\documentclass[11pt,reqno]{amsart}

\usepackage[T1]{fontenc}
\usepackage[utf8]{inputenc}
\usepackage[english]{babel}
\usepackage[letterpaper,margin=1in]{geometry}
\usepackage{amsmath}
\usepackage{amssymb}
\usepackage{amsthm}
\usepackage{graphicx}
\usepackage{booktabs}
\usepackage{float}
\usepackage{placeins}
\usepackage[ruled,vlined]{algorithm2e}
\usepackage{cite}
\usepackage[colorlinks=true,linkcolor=blue,citecolor=red,urlcolor=blue]{hyperref}

\allowdisplaybreaks
\theoremstyle{plain}
\newtheorem{theorem}{Theorem}[section]
\newtheorem{corollary}[theorem]{Corollary}
\newtheorem{lemma}[theorem]{Lemma}
\newtheorem{proposition}[theorem]{Proposition}

\theoremstyle{definition}
\newtheorem{definition}[theorem]{Definition}
\newtheorem{remark}[theorem]{Remark}

\numberwithin{equation}{section}

\newcommand{\jnt}{\displaystyle\int}
\newcommand{\jjntQ}{\displaystyle\int\!\!\!\!\int_Q}
\newcommand{\jjntomT}{\displaystyle\int\!\!\!\!\int_{]0,T[ \times \omega}}
\newcommand{\jjntomsT}{\displaystyle\int\!\!\!\!\int_{]0,T[ \times \omega_0}}
\newcommand{\Om}{\Omega}

\DeclareMathOperator{\tr}{tr}
\DeclareMathOperator{\spt}{Supp}

\title[Local null controllability of a quasi-linear system]
      {Local null controllability of a quasi-linear system\\
       and related numerical experiments}

\author[E. Fern\'andez-Cara]{Enrique Fern\'andez-Cara}
\address[E. Fern\'andez-Cara]{EDAN and IMUS, Universidad de Sevilla, Spain
  \newline\indent
  UFPB, Jo\~ao Pessoa, Brazil}
\email{cara@us.es}

\author[J. L\'imaco]{Juan L\'imaco}
\address[J. L\'imaco]{Instituto de Matem\'atica e Estat\'istica, Universidade
  Federal Fluminense, Niter\'oi, RJ, Brazil}
\email{jlimaco@id.uff.br}

\author[Y. Thamsten]{Yuri Thamsten}
\address[Y. Thamsten]{Instituto de Matem\'atica e Estat\'istica, Universidade
  Federal Fluminense, Niter\'oi, RJ, Brazil}
\email{ythamsten@id.uff.br}

\author[D. Menezes]{Denilson Menezes}
\address[D. Menezes]{Departamento de Ci\^encia e Tecnologia, Universidade
  Federal Rural do Semi-\'Arido, Cara\'ubas, RN, Brazil}
\email{denilson.menezes@ufersa.edu.br}



\date{}

\thanks{Published version: \emph{ESAIM: Control, Optimisation and Calculus of
  Variations} \textbf{29} (2023), Paper No.~27,
  \href{https://doi.org/10.1051/cocv/2023009}{doi:10.1051/cocv/2023009}.
  Received February 28, 2021; accepted January 19, 2023.
  Corresponding author: Enrique Fern\'andez-Cara (cara@us.es).}

\begin{document}

\begin{abstract}
This paper concerns the null control of quasi-linear parabolic systems where the diffusion coefficient depends on the gradient of the state variable. In our main theoretical result, with some assumptions on the regularity and growth of the diffusion coefficient and regular initial data, we prove that local null controllability holds. To this purpose, we consider the null controllability problem for the linearized system, we deduce new estimates on the control and the state and, then, we apply a Local Inversion Theorem. We also formulate an iterative algorithm of the quasi-Newton kind for the computation of a null control and an associated state. We apply this method to some numerical approximations of the problem and illustrate the results with several experiments.
\end{abstract}

\maketitle


\section{Introduction and main results} \label{intro}

   Let $d$ be an integer with $1 \leq d \leq 3$ and let $\Om \subset \mathbb{R}^d$ be a non-empty open bounded connected set, with smooth boundary $\partial \Om$.
   For a finite time horizon $T>0$, we introduce $Q := ]0,T[\times \Om$ and~$\Sigma := ]0,T[\times \partial \Om$.
   We will deal with systems where the control acts on a set of the form $]0,T[ \times \omega$, where $\omega \subset\subset \Om$ is another non-empty open set.
   More precisely, the following nonlinear parabolic system will be considered:
   \begin{equation} \label{sistema1fase}
\begin{cases}
   y_t - \nabla \cdot \left[ a(\nabla y) \nabla y \right] = \chi_\omega v \ \text{ in } Q, \\
   y = 0 \ \text{ on } \Sigma, \\
   y|_{t=0} = y_0 \ \text{ in } \Om.
\end{cases}
   \end{equation}
Here, $\chi_\omega = \chi_\omega(x)$ stands for the characteristic function of the set $\omega$ and $a : \mathbb{R} \mapsto \mathbb{R}^d$ is given.
   
   Let us provide a short physical derivation of~(\ref{sistema1fase}).
   Assuming that $d=3$ and~$\Om$ is a region occupied by a solid body, we can interpret the function $y = y(t,x)$ as the associated temperature distribution at points~$x \in \Om$ and times~$t \in \left[0,T\right]$.
   The effect of an external source acting on $\omega$ during~$[0,T]$ appears in the forcing term $\chi_\omega v$.
   Accordingly, if the mass density of the solid is assumed to be constant (for simplicity), we must have
   \begin{equation} \label{Physics}
   y_t + \nabla \cdot J = \chi_\omega v \ \text{ in } Q,
   \end{equation}
where $J$ is the heat flux density.
   In view of {\it Fourier's law,} see for instance~\cite{crank1975}, we assume that $J$ is proportional to a vector opposite to the temperature gradient, that is, there must exist a scalar $a$ (the thermal conductivity of the body) such that
   $$
J = -a\nabla y .
   $$
   If we assume that $a$ is a constant, \eqref{Physics} becomes the classical heat equation.
   For more general situations, it is pertinent to stipulate that $a = a(\nabla y)$ and this leads to~(\ref{sistema1fase}).
   
   Since the seminal work of Nash~\cite{nash1958continuity} and the velocity
   ideas introduced by De Giorgi and Moser~\cite{de1957sulla, moser1961harnack, moser1964harnack}, the theory of quasi-linear elliptic and parabolic equations has undergone an intensive development.
   Since then, using such systems to model and analyze complex physical phenomena have become usual tasks;
   see for instance~\cite{teixeira2009numerical} for a model of metallic materials with temperature-dependent parameters.
   
   
   In fluid dynamics, we also find analogous systems.
   Indeed, quasi-linear generalizations of Navier-Stokes equations are commonly used to model non-Newtonian flows, namely those with a shear-dependent viscosity;
   see~\cite{Ladyzhenskaya, malek1993non} and the recent work~\cite{berselli2019global}. In this case, the velocity field $y=y(t,x)$ and the pressure $p=p(t,x)$ of the fluid must satisfy
   $$
   y_t - \nabla\cdot\left[ \left(\nu_0 + \nu_1 |Dy|^r\right)Dy \right] + (y\cdot \nabla)y + \nabla p = f,
   $$
   for suitable constants $\nu_0,\nu_1,r > 0$, and a source $f = f(x,t)$, where $Dy$ designates the deformation tensor of $y$, that is, $Dy := \frac{1}{2}(\nabla y + \nabla y^T)$.
   
   It has also become usual to model chemotactic phenomena in population dynamics with the help of two-phase quasi-linear parabolic systems since the celebrated work~\cite{keller1970initiation}, by Keller and Segel;
   see also~\cite{hillen2009user, wang2010chemotaxis}.
   Denoting by~$u$ and~$z$ respectively the cell density and the chemical attractant concentration, the classical Keller-Segel model postulates that
   $$
   \begin{cases}
   u_t = \nabla\cdot\left( D_1 \nabla u + \chi u \nabla z \right) + k,\\
   z_t = D_2\Delta z + g - h.
   \end{cases}
   $$
   Here, $D_1,D_2>0$ are diffusion coefficients, $k$ models the birth/death rate of cells and~$g$ (resp.~$h$) represents the production rate (resp.~the degradation rate) of the chemical substance.
   
   We can also mention systems of the kind~\eqref{sistema1fase} in image processing modelling.
   For instance, the so called Perona-Malik (anisotropic diffusion) model~\cite{perona1990scale} reads
   $$
   I_t - \nabla\cdot \left[ a(\nabla I) \nabla I \right] = 0 \ \text{ in } \ Q,
   $$
   where $I(t,\cdot)$ represents a family of gray scale images over the planar set $\Om$ and either
   $$
   a(p) = e^{-|p|^2/K^2},
   $$
   or
   $$
   a(p) = 1/[1+(|p|/K)^2]
   $$
for some constant $K>0$.
   
   In the last decades, there have been many advances in the control of linear and semilinear parabolic systems;
   see for instance~\cite{fabre1995approximate, fernandez2000null, Fursikov-Imanuvilov, lebeauR} and the references therein.
   The works~\cite{bonifacius2018second, casas2018analysis, casas1995optimal, pan2001optimal} dealt with the optimal control of certain quasi-linear parabolic systems.
   Recent efforts concerning the internal controllability of systems with diffusion coefficients with nonlocal dependence on the state variable and/or its gradient have appeared in~\cite{CCLM, fernandez2015theoretical};
   for chemotaxis models, we also have~\cite{chaves2015uniform, chaves2017controllability};
   some exact controllability properties for nonlinear parabolic equations of quasi-linear kind can be found in~\cite{fernandez2021theoretical, CHCV, LiuZhang}.
   Note that, for many of these models, internal controllability results lead to boundary control after a relatively simple argument relying on domain modification.
   
   As indicated in~\cite{LiuZhang}, in order to address the case where the coefficient also depends on the gradient, specific techniques are needed.
   However, as the previous discussion shows, this is a very important class of models from a practical perspective.
   Therefore, the control analysis of systems like~\eqref{sistema1fase} is of great interest.

   In the present work, we extend the methods introduced in~\cite{fernandez2021theoretical} and we accept diffusion coefficients that depend on the state gradient.
   The idea is, as in other works, to first consider the linearized at zero of~\eqref{sistema1fase} and prove a related null control result.
   Then, this is used in the context of a local inversion theorem to get the desired (local) controllability property.
   Nevertheless, in view of the complexity of the equation, to make this argument work, we have to establish a collection of nontrivial estimates on the solutions to the linear control problems.
   To this respect, we apply a boot-strapping argument:
   after requiring higher regularity to the initial data, the controls can be found with higher regularity, thus guaranteeing more regularity to the state variables;
   in turn, this makes it possible to go further and get an even more regular control, a more regular state, etc.
   The authors of~\cite{duprez2019bilinear} employ a similar idea.
  
  The need of ensuring sufficiently regular state and control properties forces to constrain the spatial dimension $d$ to be at most three;
   to this regard, see~Remark~\ref{rem3.6}.
   
   Throughout the whole paper, the following hypotheses will be assumed:
   
\begin{itemize}
   \item[\textbf{H1}] The function $a : \mathbb{R}^d \mapsto \mathbb{R}$ is of class~$C^4$ in~$\mathbb{R}^d$.
   \item[\textbf{H2}] There exist constants $K_0$, $a_0 > 0$ and~$r\geq 1$, such that $a$ satisfies
   $$
   \begin{cases}
   a(q) \geq a_0 > 0, \\ 
   |D^\beta a(q)| \leq K_0\left(1+|q|^{(r-|\beta|)^+}\right),
   \end{cases}
   $$
   for all $q \in \mathbb{R}^d $ and any multi-index $\beta$ with~$|\beta| \leq 4$.
   Here and in the sequel,  $z^+$ stands for the positive part of the real number~$z$, that is, $z^+ := \max\left(z, 0\right)$.
\end{itemize}

   An example of a function for which both~\textbf{H1} and~\textbf{H2} hold is $a(p) = a_0 + a_1 |p|^r$, for positive constants~$a_0$ and~$a_1$ and~sufficiently large $r \geq 1$.
   
   Let us specify the controllability property considered in this paper:
   
\begin{definition} \label{definitionNC}
   It is said that $(\ref{sistema1fase})$ is locally null-controllable at time $T$ if there exists $\eta > 0$ such that, if the initial data satisfies
   \begin{equation} \label{1.2p}
y_0 \in H^5(\Om) \cap H^1_0(\Om)
   \end{equation}
and the compatibility conditions
\begin{equation}
    \Delta y_0,\Delta^2 y_0 \in H^1_0(\Om)
\end{equation}
and is sufficiently small in the sense that
\begin{equation}\label{1.2pp}
    \|y_0\|_{H^5} \leq \eta,
\end{equation}
there exist controls $v \in L^2(]0,T[\times \omega)$ and associated solutions to~$(\ref{sistema1fase})$ with
   \begin{equation} \label{NC}
   y(T,x) = 0 \ \text{ in } \ \Om.
   \end{equation}
\end{definition}

\begin{remark}\label{Novelty}
   To the best of our knowledge, \cite{de2019local} and the present paper are the first works dealing with the null controllability of quasi-linear parabolic PDE's where the diffusion coefficient can be, among other possibilities, a general power law.
\end{remark}

   Our main result is the following:

\begin{theorem}\label{Maintheorem}
   For every $T>0$, $(\ref{sistema1fase})$ is locally null-controllable at time $T$.
\end{theorem}

   Using standard techniques, it is not difficult to deduce from this theorem a similar boundary controllability result:
   
\begin{corollary}\label{corol-1.4}
   For any $T>0$, there exists $\eta > 0$ such that, for any $y_0$ satisfying~\eqref{1.2p}--\eqref{1.2pp}, there exist controls $h \in L^2(\Sigma)$ and associated solutions to
   \begin{equation} \label{bdryControl}
   \begin{cases}
   y_t - \nabla \cdot \left[ a(\nabla y) \nabla y \right] = 0 \ \text{ in } Q, \\
   y = h \ \text{ on } \Sigma, \\
   y|_{t=0} = y_0 \ \text{ in } \Om
   \end{cases}
   \end{equation}
such that~\eqref{NC} holds.
   Furthermore, given a non-empty open set $\gamma \subset \partial\Om$, it is possible to select $h$ such that $\spt(h) \subseteq \left[0,T\right]\times \gamma$.
\end{corollary}

   For the proof of Theorem~\ref{Maintheorem}, we will rewrite the null controllability problem in the form
   \begin{equation} \label{equation_H}
H(y,v) = (0,y_0), \quad (y,v) \in Y,
   \end{equation}
where the Hilbert space $Y$ of state-control pairs $(y,v)$ will be chosen appropriately
   (in particular, to any couple $(y,v) \in Y$ we will require to satisfy $y(T\,,\cdot) = 0$) and~$H : Y \mapsto Z$ will be an adequate nonlinear mapping ($Z$ is another Hilbert space).
   Among other things, we prove that $H$ is well-defined and strictly differentiable at $(0,0)$, with a surjective derivative $H'(0,0) \in \mathcal{L}(Y;Z)$.
   Accordingly, it will be possible to apply a local inversion theorem to~\eqref{equation_H} and deduce that, for any ``small'' $y_0$, we can steer \eqref{sistema1fase} to zero.

   Our original contribution in this context is to provide good definitions of the spaces~$Y$ and~$Z$, as well as of the mapping~$H$.
   We must harmonize conflicting relations when specifying~$Y$ and~$Z$:
   indeed, they must be sufficiently small to guarantee the correct definition and needed regularity of~$H$ and, also, sufficiently large to ensure that the solutions to the linearized problem belong to~$Y$.
   Furthermore, we can only establish the well-posedness and differentiability of~$H$ after several new weighted energy estimates for the solutions to the linearized system.
   
   Note that it seems difficult to assert the global null controllability of~\eqref{sistema1fase}.
   We could try to apply a global inversion argument, like for instance in~\cite{LasTrig}, based on the existence and boundedness of~$H'(\bar{y},\bar{v})$ for all~$(\bar{y},\bar{v})$.
   However, this requires in practice the resolution of null controllability problems for systems of the form
   \[
\begin{cases}
   y_t - \nabla \cdot \left[ a(\nabla \bar{y}) \nabla y + (a'(\nabla \bar{y})\cdot\nabla y) \nabla\bar{y} \right] = \chi_\omega v \ \text{ in } Q, \\
   y = 0 \ \text{ on } \Sigma, \\
   y|_{t=0} = y_0 \ \text{ in } \Om.
\end{cases}
   \]
with good estimates of~$v$ and~$y$, uniform with respect to~$\bar y$.
   And this is an open and probably difficult problem.
   
   We could also try to apply a fixed-point argument, that is, search for a fixed-point of the mapping $\tilde y \mapsto y$, where $y$ is now, together with some~$v$, a solution to the null controllability problem for
   \[
\begin{cases}
   y_t - \nabla \cdot \left[ a(\nabla \tilde{y}) \nabla y \right] = \chi_\omega v \ \text{ in } Q, \\
   y = 0 \ \text{ on } \Sigma, \\
   y|_{t=0} = y_0 \ \text{ in } \Om,
\end{cases}
   \]
assuming that this mapping can be well defined.
   But, unfortunately, it is difficult to find (and not clear at all) a space where this can be done.

   In the second part of the paper, we will be concerned with the computation of null controls for~\eqref{sistema1fase}.
   In the literature on the subject, the numerical controllability of linear and nonlinear PDE's has been considered in many papers;
   for instance, see~\cite{B,BHL,GL,fernandez2014numerical,LT,munch2010numerical} and the references therein.
   Here, we will argue as in~\cite{CCLM, CHCV}, taking advantage of the surjectivity of $H'(0,0)$.

   Thus, let~$Y$ be the Hilbert space where we can find a solution $(y,v)$ to~\eqref{equation_H} (see~\eqref{spaceY} for details).
   We introduce the following iterative algorithm, in the sequel denoted \textbf{ALG~1}:

\begin{algorithm}[H] \label{ALG1}
\SetAlgoLined
Start with $n=0$, an initial guess $(y^0,v^0) \in Y$, the current error $\epsilon$ and tolerance $\epsilon_0 > 0$.

\While{$\epsilon\geq \epsilon_0$}{
1: Compute 
   $$
(y^{n+1}, v^{n+1}) =(y^n, v^n) - H'(0,0)^{-1} (H(y^n, v^n) - (0,y_0));
   $$ \\
2: Update $\epsilon$:
   $$
\epsilon \gets \| (y^{n+1},v^{n+1}) - (y^n,v^n) \|_{L^2(Q)}/\|(y^n,v^n)\|_{L^2(Q)};
   $$\\
3: Update $n$:
   $$
n \gets n+1.
   $$\\
}
\Return{$(y^{n+1},v^{n+1})$}
\caption{Numerical algorithm for computing a control, and the corresponding state, for \eqref{sistema1fase}}
\end{algorithm}


   In these iterates, we use~$H'(0,0)^{-1}$, which is by definition an inverse to the left of~$H'(0,0)$.
   The precise definition will be given in Section~\ref{NL-localNC}.
   
   We remark that \textbf{ALG~1} is an elementary quasi-Newton method;
   below, it will be shown to converge in the space~$Y$.
   We will analyze the convergence and the computational cost of \textbf{ALG~1} in~Section~\ref{Sec4} and we will present several numerical experiments to illustrate our results.

%

   The plan of the paper is the following:

\begin{itemize}

   \item \textbf{Study of the linearized system.}
   In Section~\ref{SecNCforLinear1Phase}, we carry out a detailed study of the linearization of~(\ref{sistema1fase}) around zero.
   We recall Fursikov-Imanuvilov's method and solve the corresponding null controllability problem.
   We employ techniques that make it possible to find controls with high regularity, as long as the initial state is sufficiently smooth.
   
   \item \textbf{Application of a Local Inversion Theorem.}
   In Section~\ref{SecLocNC1Phase}, with the help of an adequate functional-analytic setting, the controllability problem is written in the form~\eqref{equation_H}.
   At this point, we need higher order estimates for the solutions to non-homogeneous linearized problems;
   as already mentioned, this is one of the main contributions of the paper.
   Then, we find that the relevant mapping is well-defined and differentiable (in an appropriate sense), whence we can apply a local inversion result, taken from~\cite{OptimalControl}.
   This yields the desired result.
   
   \item \textbf{Computation of null controls.}
   We present in~Section~\ref{Sec_Comments} some numerical methods for the computation of a solution to the null controllability problem for~(\ref{sistema1fase}).
   The results of some numerical experiments in spatial dimensions $d=1$ and~$d=2$ are also given.
   
   \item \textbf{Comments and concluding remarks.}
   We summarize the achievements in Section~\ref{Sec_Comments}.
   There, we have also included several additional comments and pertinent open questions.
   
\end{itemize}

   In principle, throughout this paper, we regard all derivatives in the distributional sense.
   For any~$d$-dimensional multi-index $\beta$, the symbol $D^\beta$ stands for the corresponding spatial differential operator;
   explicitly,
   $$
D^\beta = \partial_1^{\beta_1}\cdots \partial_d^{\beta_d}.
   $$
   
   We will use $\|\cdot\|$ to denote the standard norm in~$L^2(\Om)$.
   In general, for any other Banach space~$E$, the corresponding norm will be denoted by $\|\cdot \|_E$.
   The symbol $C$ will stand for a generic positive constant, depending on~$\Om$, $\omega$, $T$ and the other data appearing in hypotheses~\textbf{H1} and~\textbf{H2}.
   For any given nonnegative integer $k$ and any function $g$, whenever it makes sense, we set
   $$
|D^k g| := \left( \sum_{|\beta| \leq k} |D^\beta g|^2\right)^{1/2}.
   $$
   Finally, the integration elements $dt$, $dx$, etc.\ will be usually omitted.

\section{Study of the linearized problem} \label{SecNCforLinear1Phase}

   The linear system considered in this section is the following:
   \begin{equation} \label{Sistema1FaseLinear}  
\begin{cases}
   y_t - \sigma\Delta y  = \chi_\omega v + f\ \text{ in } Q, \\
   y = 0\ \text{ on } \Sigma, \\
   y|_{t=0} = y_0\ \text{ in } \Om,
\end{cases}
   \end{equation}
where $\sigma$ is a positive constant.
   We must also consider the corresponding adjoint system, given by
   \begin{equation} \label{Sistema1FaseLinearAdj}
\begin{cases}
-\varphi_t - \sigma\Delta \varphi = F\ \text{ in } Q, \\
\varphi = 0\ \text{ on } \Sigma, \\
\varphi|_{t=T} = \varphi^T\ \text{ in } \Om.
\end{cases}
   \end{equation}
   It is well-known that the null controllability of (\ref{Sistema1FaseLinear}), together with an estimate of the control, is equivalent to the observability of (\ref{Sistema1FaseLinearAdj}).

   Let $\omega_0 \subset\subset \omega$ be a non-empty open set.
   The first result we recall is due to Fursikov and Imanuvilov (cf.~\cite{Fursikov-Imanuvilov}).
   It is fundamental in this work:
	
\begin{lemma} \label{fursikovslemma}
   There exists a function $\alpha_0 \in C^4(\overline{\Om})$ such that
   $$
\alpha_0 > 0 \ \text{ in } \ \Om, \ \ \alpha_0 \equiv 0 \ \text{ on } \ \partial \Om \ \text{ and } \  |\nabla \alpha_0|> 0
\ \text{ in } \ \overline{\Om} \backslash \omega_0.
   $$
\end{lemma}

   Let us fix $\ell \in C^\infty\left(\left[0,T\right]\right)$ with
   $$
\ell(t) \geq \max\left(T^2/8,\ t(T-t) \right) \ \text{ in } \left[0,T/2\right]
\ \text{ and } \ \ell(t) = t(T-t) \ \text{ in } \ \left[T/2,T\right]
   $$
and let us set
   $$
\zeta(t,x) := \frac{e^{\lambda \alpha_0(x)}}{\ell(t)}
\ \text{ and } \ \alpha(t,x) := \frac{e^{R\lambda}-e^{\lambda\alpha_0(x)}}{\ell(t)} = \frac{\overline{\alpha}(x)}{\ell(t)} ,
   $$
where $R > \|\alpha_0\|_{L^\infty} + \log(2)$ and~$\lambda > 0$.
   We will need in the sequel a global Carleman inequality, where we estimate the quantity
   $$
I(s,\lambda;\xi) := \jjntQ e^{-2s\alpha}\left[(s\zeta)^{-1}\left(|\xi_t|^2 + |\Delta \xi|^2 \right) + \lambda^2(s\zeta)|\nabla \xi|^2 + \lambda^4(s\zeta)^3|\xi|^2 \right] .
   $$
   The result is the following:
   
\begin{proposition} \label{Carleman}
   There exist positive constants $\lambda_0$, $s_0$ and~$C_0$, depending only on $\Om$, $\omega_0$, $\sigma$ and~$T$, such that for any $s\geq s_0$, $\lambda \geq \lambda_0$, $F\in L^2(Q)$ and~$\varphi^T \in L^2(\Om)$, the solution to~$(\ref{Sistema1FaseLinearAdj})$ corresponding to the data $F$ and~$\varphi^T$ satisfies
   \begin{equation} \label{CI-new}
   I(s,\lambda;\varphi) \leq C_0\left(\jjntQ e^{-2s\alpha}|F|^2\  + \jjntomsT e^{-2s\alpha}(s\zeta)^3|\varphi|^2 \right).
   \end{equation}   
\end{proposition}

   The proof is given in~\cite{CCLM}.
   From now on, we fix $s = s_0$ and $\lambda = \lambda_0$.

   We will now prove that the linearized system~\eqref{Sistema1FaseLinear} is null-controllable, with regular controls.
   Let us set 
   $$
\alpha_1 := \min_{x \in \overline{\Om}}\overline{\alpha}(x) \ \text{ and } \ \alpha_2 := \max_{x \in \overline{\Om}}\overline{\alpha}(x). 
   $$
   We will assume that $R$ is sufficiently large to have $2\alpha_1 \geq \alpha_2$.
   We note that this also implies $(r+1)\alpha_1 \geq \alpha_2$, since $r \geq 1$ (cf. \textbf{H2} in Section \ref{intro}).
   Thus, we have:
   \begin{equation} \label{eq1}
e^{s\alpha_1(x)/\ell(t)} \leq  e^{s\alpha(x)/\ell(t)} \leq  e^{2s\alpha_2(x)/\ell(t)} \leq  e^{2s\alpha_1(x)/\ell(t)} .
   \end{equation}
   The following weights will be needed: 
   \begin{equation}\label{def-rho}
\rho_k(t,x) := e^{s\alpha(x)} \ell(t)^{k/2} \hspace{1.0cm} (k=0,1,2,...); \ \text{ in particular, we set $\rho:=\rho_0$.}
   \end{equation}

\begin{theorem} \label{ControlOfLinearProblem}
   Let us assume that $\rho_3 f \in L^2(Q)$ and $y_0 \in L^2(\Om)$.
   There exists a control $v \in L^2\left(]0,T[\times \omega\right)$ which, together with the solution to~$(\ref{Sistema1FaseLinear})$ corresponding to~$v$, $f$ and~$y_0$, satisfies:
   \begin{equation} \label{EstLinControl1}
\jjntQ\rho^2 |y|^2 + \jjntomT \rho_3^2 |v|^2 \leq C\left( \|y_0\|^2 + \jjntQ \rho_3^2 |f|^2  \right).
   \end{equation}
   In particular, $v$ is a control that drives the solution to~\eqref{Sistema1FaseLinear} from~$y_0$ exactly to zero at time~$T$.
   Furthermore, we can choose $v$ satisfying 
   $$
(\rho_7 v)_t \in L^2\left(]0,T[\times \omega\right)
\ \text{ and } \ \rho_7 v \in L^2\left(0,T;H^2(\omega)\cap H^1_0(\omega)\right),
   $$
together with the estimates
   \begin{equation} \label{EstLinControl2}
\jjntomT \left[ |(\rho_7 v)_t|^2+ \left|\Delta(\rho_7 v)\right|^2\right]
\leq C\left(\|y_0\|^2 + \jjntQ \rho_3^2 |f|^2  \right).
   \end{equation}
\end{theorem}

   This result is well-known.
   The proof of the first part is due to Fursikov and~Imanuvilov~\cite{Fursikov-Imanuvilov}.
   For the proof of~\eqref{EstLinControl2}, see for instance~\cite{fernandez2021theoretical}.
   
   For future purposes, let us briefly recall the way the null control $v$ is found in~\cite{Fursikov-Imanuvilov}.
   Thus, let $\chi \in C^\infty_c(\omega)$ be given, with~$0 \leq \chi \leq 1$ and~$\chi|_{\omega_0} \equiv 1$ and let us introduce the linear space 
   $$
 P_0 := \left\{ w \in C^2(\overline{Q}) : w = 0 \ \text{ on } \Sigma \right\}
   $$
 and the bilinear form
   $$
 \pi(w,\widetilde{w}) := \jjntQ \rho^{-2}L^*w L^*\widetilde{w}\  + \jjntQ \rho_3^{-2}\chi w\widetilde{w} ,
   $$
 where $L^* w := -w_t - \sigma\Delta w$.
   In view of the {\it unique continuation property} of the heat operator, it is clear that~$\pi(\cdot\,,\cdot)$ is a scalar product in~$P_0$.
   Let~$P$ be the completion of $P_0$ for the inner product $\pi(\cdot\,,\cdot)$.
   From the Carleman inequality (i.e. Proposition \ref{Carleman}), together with the energy estimates for the solution to~(\ref{Sistema1FaseLinearAdj}) (which can be deduced in a standard way), we see that the linear form
   $$
 m(w) := \jnt_\Om y_0(x) w(0,x)\,dx + \jjntQ fw
   $$
 is continuous on~$P$.
   Consequently, from the {\it Riesz Representation Theorem,} it follows that there exists a unique $p$ satisfying
   \begin{equation} \label{2.9p}
 \pi(p,w) = m(w) \quad \forall w \in P, \quad p \in P.
   \end{equation}
   
   Now, we take
   $$
 y := \rho^{-2}L^*p, \ \ \hat{v}:= -\rho_3^{-2}p \ \text{ and } \ v := \bigl. \chi\hat{v} \bigr|_{]0,T[ \times \omega}
   $$
and we find that~$(y,v)$ is a state-control couple of (\ref{Sistema1FaseLinear}), with $y(0\,,\cdot) = y_0$.
   Furthermore, we also have
   $$
 \pi(p,p) \leq C\left( \|y_0\|^2 + \jjntQ \rho_3^2 |f|^2  \right)^{1/2}\pi(p,p)^{1/2} ,
   $$
 whence we get~(\ref{EstLinControl1}).

   Note that, in view of~\eqref{EstLinControl1}, the control $v$ also fulfills the following regularity properties:
   $$
(\rho_7 v)_t \in L^2\left(]0,T[\times \omega\right)
\ \text{ and } \ \rho_7 v \in L^2\left(0,T;H^2(\omega)\cap H^1_0(\omega)\right).
   $$

\begin{corollary} \label{RegularityForLinear1PhaseProblem}
   Let the assumptions in~Theorem~\ref{ControlOfLinearProblem} be satisfied and let~$v$ be a null control for~\eqref{Sistema1FaseLinear} furnished by~Theorem~\ref{ControlOfLinearProblem} satisfying~\eqref{EstLinControl2}.
   The following holds:
   
\begin{itemize}

   \item[(a)] The associated state $y$ satisfies
   \begin{equation} \label{regularity1}
   \begin{array}{l} \displaystyle
         \sup_{\left[0,T\right]} \left[\jnt_\Om \rho_5^2 |y|^2 \,dx \right](t) + \jjntQ \rho_5^2 \left|\nabla y\right|^2  
	   \leq C\bigg(\|y_0\|^2 + \jjntomT \rho_3^2 |v|^2  + \jjntQ \rho_3^2 |f|^2  \bigg). 
   \end{array}   
   \end{equation}
   
   \item[(b)] If we assume that $y_0 \in H^1_0(\Om)$, then we also have
   \begin{equation} \label{regularity2}
   \begin{array}{l} \displaystyle
      \sup_{\left[0,T\right]} \left[ \jnt_\Om \rho_{7}^2 |\nabla y|^2 \,dx \right] \!(t) 
      \!+\!\! \jjntQ \!\rho_{7}^2 \!\left[ |y_t|^2 \!+\! |\Delta y|^2 \right] 
	\leq C\!\left( \! \|y_0\|_{H^1}^2 \!+\!\!\jjntomT \!\rho_3^2|v|^2  \!+\! \jjntQ \rho_3^2 |f|^2 \!\right).
	\end{array}   
   \end{equation}
   Moreover, we have the following regularity properties in this case:
   \[
(\rho_{15} v_t)_t \in L^2(]0,T[\times \omega), \ \ \rho_{15} v_t \in L^2(0,T;H^2(\omega)), \ \ \rho_{15}v \in L^2(0,T;H^4(\omega))
   \]
and
   \begin{equation} \label{regularity3}
   \begin{array}{l} \displaystyle
   \jjntomT \left[|(\rho_{15}v_t)_t|^2 + |\Delta(\rho_{15}v_t)|^2 \right]  + \int_0^T \|\rho_{15}v\|_{H^4(\Om)}^2
   \leq C\left( \|\nabla y_0\|_{H^1}^2 + \jjntQ \rho_3^2 |f|^2 \right).
   \end{array}
   \end{equation}
   
   \item[(c)] Furthermore, if $y_0 \in H^3(\Om)\cap H^1_0(\Om)$, $\Delta y_0 \in H^1_0(\Om)$, $\rho_9 f_t \in L^2(Q)$ and~$f(0\,,\cdot) \in H^1_0(\Om)$, then we have:
   \begin{equation} \label{regularity4}
   \begin{array}{l} \displaystyle
      \sup_{\left[0,T\right]} \left[ \jnt_\Om \left( \rho_9^2 |y_t|^2 
      + \rho_{11}^2 |\Delta y|^2 + \rho_{11}^2 |\nabla y_t|^2 \right)\,dx \right](t)
	   + \jjntQ \left( \rho_9^2 |\nabla y_t|^2  + \rho_{11}^2 |y_{tt}|^2  +\rho_{11}^2 |\Delta y_t|^2 \right)
	   \\ \noalign{\smallskip} \displaystyle
	   \hspace{0.2cm} \ \leq C\bigg\{ \|y_0\|_{H^3}^2 
	   \!+\! \jjntomT \left[ \rho_3^2 |v|^2 \!+\! |(\rho_7v)_t|^2 \right]
	   \!+\! \|f(0\,,\cdot)\|_{H^1}^2 \!+\! \jjntQ\left( \rho_3^2 |f|^2 \!+\! \rho_9^2 |h_t|^2\right)  
	   \bigg\}.
   \end{array}   
   \end{equation}
   
   \item[(d)] Finally, under the assumptions in part~(c), whenever $y_0 \in H^5(\Om)\cap H^1_0(\Om)$, $\Delta y_0 \in H^1_0(\Om)$, $\Delta^2 y_0 \in H^1_0(\Om)$, $\rho_7 f \in L^2(0,T; H^2(\Om))$, $\rho_{11} f_t \in L^2(0,T; H^1_0(\Om))$, $\rho_{17}f_{tt} \in L^2(Q)$, $f = 0$ on~$\Sigma$, $f(0\,,\cdot) \in H^3(\Om)\cap H^1_0(\Om)$, $\Delta f(0) \in H^1_0(\Om)$ and~$f_t(0\,,\cdot) \in H^1_0(\Om)$, one has:
   \begin{equation} \label{regularity5}
   \begin{array}{l} \displaystyle
      \sup_{\left[0,T\right]} \left[\jnt_\Om \left( \rho_{15}^2 |y_{tt}|^2 
      + \rho_{17}^2|\nabla y_{tt}|^2 + \rho_{17}^2(\Delta y_t)^2 \right) \,dx \right](t)
      \\ \noalign{\smallskip} \displaystyle
	   \hspace{.5cm} + \jjntQ \left[ \rho_{15}^2 |\nabla y_{tt}|^2  
	   + \rho_{17}^2 |y_{ttt}|^2  + \rho_{17}^2 |\Delta y_{tt}|^2   \right]
	   + \jnt_0^T \left[ \left\|\left( \rho_{11}  y\right)(t,\,\cdot) \right\|_{H^4(\Om)}^2 + \|(\rho_{13} y_t)(t,\,\cdot)\|_{H^3}^2 \right]
	   \\ \noalign{\smallskip} \displaystyle
	   \hspace{0.1cm}
	   \leq C\Big\{ \|y_0\|_{H^5}^2 \!+\!\! \jjntomT \!\! [\rho_3^2 |v|^2 \!+\! |\Delta(\rho_7v)|^2 \!+\! |(\rho_7 v)_t|^2
	   \!+\! |(\rho_{15}v_t)_t|^2 \!+\! |\Delta(\rho_{15}v_t)|^2 ] 
	   \\ \noalign{\smallskip} \displaystyle
	   \hspace{0.5cm}
	    \!+\! \int_0^T \|\rho_{15}v(t,\,\cdot)\|_{H^4(\Om)}^2 
	   + \!\! \jjntQ [\rho_3^2 |f|^2 \!+\! \rho_7^2 |\Delta f|^2 \!+\! \rho_9^2 |f_t|^2 
	   \!+\! \rho_{11}^2 |\nabla f_t|^2 \!+\! \rho_{17}^2 |f_{tt}|^2] 
	   \\ \noalign{\smallskip} \displaystyle
	   \hspace{7.7cm}
	   + \ \|f(0\,,\cdot)\|_{H^3}^2 + \|f_t(0\,,\cdot)\|_{H^1_0}^2\Big\}.
   \end{array}   
   \end{equation}
   In particular, $\rho_{13} y \in L^\infty(0,T;H^3(\Om))$ and
   \begin{equation} \label{regularity6}
\begin{array}{l} \displaystyle
   \sup_{\left[0,T\right]}\|\rho_{13}y\|_{H^3}^2(t) \leq C\jnt_0^T\left(\|(\rho_{13} y)(t\,,\cdot)\|_{H^4}^2 
   + \|(\rho_{13}y)_t(t\,,\cdot)\|_{H^2}^2\right) \,dt
   \leq C S(y_0,v,f),
\end{array}   
   \end{equation}
where $S(y_0,v,f)$ stands for the right hand side of~\eqref{regularity5}.

\end{itemize}

\end{corollary}

   The proof of this result is given in Appendix~\ref{Sec-App-A}.
      
\section{The local null controllability of (\ref{sistema1fase})} \label{SecLocNC1Phase}

   In this section, we will rewrite the local null controllability of (\ref{sistema1fase}) as a nonlinear equation in a Hilbert space.
   Then, we will apply an ``abstract'' local inversion result.
   This will lead to the controllability result in~Theorem~\ref{Maintheorem}.
   
\subsection{A general local inversion theorem}

   In this section, $Y$ and~$Z$ will stand for generic Hilbert spaces.
   We will denote the closed ball centered at~$0$ of radius $r$ in~$Y$ (resp.~$Z$) by $B_Y(0;r)$ (resp.~$B_Z(0;r)$).

   We begin by recalling the definition of strict differentiability:
      
\begin{definition} \label{StrictlyDiff}
   A mapping $H : Y \mapsto Z$ is strictly differentiable at $0$ if there exists $\Lambda \in \mathcal{L}(Y;Z)$ such that, for any $\epsilon>0$, there exists $\delta(\epsilon) > 0$ with
   $$
\|H(y) - H(\overline{y}) - \Lambda(y-\overline{y})\|_Z \leq \epsilon\|y-\overline{y}\|_Y
\quad \forall y,\overline{y} \in B_Y(0;\delta(\epsilon)).
   $$
   If this is the case, then $\Lambda$ is unique and is usually denoted by~$H'(0)$.
\end{definition}

   The following result holds:

\begin{theorem} \label{Liusternik}
   Let $H : Y \mapsto Z$ be given with $H(0) = 0$.
   Suppose that $H$ is strictly differentiable at $0$ and $H^\prime(0)$ is onto.
   Then, there exist $\delta > 0$, $r > 0$ and~$W : B_Z(0;\delta) \mapsto B_Y(0;r)$ such that
   $$
H(W(z)) = z \quad \forall z \in B_Z(0;\delta).
   $$
   Furthermore, if there exists $M > 0$ such that
   $$
\| H'(0)y \|_Z \geq M \|y\|_Y \quad \forall y \in Y,
   $$
then we can take $\delta = \delta(\epsilon)$ and $r = (M^{-1} - \epsilon)^{-1} \delta(\epsilon)$, where $\delta(\epsilon)$ is as in~Definition~\ref{StrictlyDiff} and $\epsilon$ is any number satisfying $0 < \epsilon < M^{-1}$.
\end{theorem}

   We can find a proof of this result in~\cite{OptimalControl}.
   
   The particular statement we give here is similar to the one in~\cite{chaves2015uniform}.
   Obviously, we must view the mapping $W$ as an {\it inverse to the right} of~$H$.

\subsection{Proof of the local null controllability}\label{NL-localNC}

   Let us introduce the weighted $L^2$-spaces $H_k := L^2(Q,\mu_k)$.
   They will be of help in the definition of the particular spaces $Y$ and~$Z$ needed in the proof of~Theorem~\ref{Maintheorem}.
   Here, $\mu_k$ is the absolutely continuous measure in~$Q$ whose Radon-Nikod\'ym derivative is $\rho_k^2$, that is, with integration element $d\mu_k = \rho_k^2 \,d(t,x)$.
   Thus, a (class of) function(s) $v : Q \mapsto \mathbb{R}$ belongs to~$H_k$ if and only if it is measurable and
   \[
\jjntQ |v|^2 d\mu_k = \jjntQ \rho_k^2 |v|^2 < +\infty.
   \]
   The space $H_k$ is endowed with the natural norm
   \begin{equation}\label{3.1a}
\|u\|_k := \left(\jjntQ |u|^2 d\mu_k \right)^{1/2} = \left( \jjntQ \rho_k^2 |u|^2 \right)^{1/2} .
   \end{equation}

   We will set
   \begin{equation} \label{spaceY}
\begin{array}{l}
   Y := \{ (y,v) \in H_0 \times H_3 : y \in L^2(0,T;H^1_0(\Om)), 
   \ y_t \in L^2(Q),
   \ y(0\,,\cdot) \in H^5(\Om)\cap H^1_0(\Om),
   \\ \noalign{\smallskip} \displaystyle
   \phantom{Y := \{ } \Delta y(0,\,\cdot),\Delta^2 y(0,\,\cdot) \in H^1_0(\Om),\rho_7 v \in L^2(0,T;H^2(\omega)\cap H^1_0(\omega)),
   \\ \noalign{\smallskip} \displaystyle
   \phantom{Y := \{ } \rho_{15} v \in L^2(0,T;H^4(\omega)),
   \ \rho_7 v_t \in L^2(]0,T[\times\omega),
   \\ \noalign{\smallskip} \displaystyle
   \phantom{Y := \{ } \rho_{15}v_t \in L^2(0,T;H^2(\omega)),
   \ \rho_{15}v_{tt} \in L^2(]0,T[\times \omega) ,
   \\ \noalign{\smallskip} \displaystyle
   \phantom{Y := \{ } f:= y_t-a(0)\Delta y - \chi_\omega v \ \text{ satisfies } 
   \ f \in H_3, \ f_t \in H_9, \  f_{tt} \in H_{17},
   \\ \noalign{\smallskip} \displaystyle
   \phantom{Y := \{ } \ \rho_7 f \in L^2(0,T;H^2(\Om)),
   \ \rho_{11} f_t \in L^2(0,T;H^1_0(\Om)),
   \\ \noalign{\smallskip} \displaystyle
   \phantom{Y := \{ } f = 0 \ \text{ on } \Sigma,
   \ f(0\,,\cdot) \in H^3(\Om) \cap H^1_0(\Om),
   \ \Delta f(0\,,\cdot) \in H^1_0(\Om),
   \ f_t(0\,,\cdot) \in H^1_0(\Om) \}.
\end{array}.
   \end{equation}
   
   This is a Hilbert space for the norm
   \begin{equation} \label{NormOfY}
\begin{array}{l}
   \|(y,v)\|_Y^2 := \|v\|_3^2 + \|v_t\|_7^2 + \|\Delta v\|_7^2 + \|D^4 v\|_{15}^2 + \|\Delta v_t\|_{15}^2 + \|v_{tt}\|_{15}^2 
   \\ \noalign{\smallskip} \displaystyle
   \hspace{1.9cm}\ +\ \|f\|_3^2  + \|f_t\|_9^2 + \|f_{tt}\|_{17}^2 + \|\Delta f\|_7^2 + \|\nabla f_t\|_{11}^2 
   \\ \noalign{\smallskip} \displaystyle
   \hspace{1.9cm}\ +\ \|f(0\,,\cdot)\|_{H^3}^2 + \|f_t(0\,,\cdot)\|_{H^1_0}^2 
   + \|y(0\,,\cdot)\|_{H^5}^2
\end{array}   
   \end{equation}
where, as before, we have set $f := y_t - a(0)\Delta y - \chi_\omega v$ and we take the norms of $v$ and its derivatives over~$]0,T[\times\omega$.
   In order to define the space $Z$, let us introduce
   \begin{equation} \label{spaceZ}
\begin{array}{l}
   F := \{ g \in H_3 : \rho_7 g \in L^2(0,T;H^2(\Om)),
   \ g_t \in H_9,
   \ \rho_{11} g_t \in L^2(0,T;H^1(\Om)),
   \ g_{tt} \in H_{17},
   \\ \noalign{\smallskip} \displaystyle
   \phantom{F := \{ } g = 0 \ \text{ on } \Sigma,
   \ g(0\,,\cdot) \in H^3(\Om)\cap H^1_0(\Om),
   \ g_{t}(0\,,\cdot) \in H^1_0(\Om) \}
\end{array}
   \end{equation}
and the Hilbertian norm
   \begin{equation} \label{NormOfF}
\begin{array}{l}
   \|g\|_F^2 := \|g\|_3^2 + \|\Delta g\|_7^2 + \|g_t\|_9^2 + \|\nabla g_t\|_{11}^2 + \|g_{tt}\|_{17}^2
   + \|g(0\,,\cdot)\|_{H^3}^2 + \|g_t(0\,,\cdot)\|_{H^1}^2 .
\end{array}
   \end{equation}
   We will set
   \begin{equation}\label{3.1aa}
Z := F\times\left[H^5(\Om)\cap H^1_0(\Om)\right],
   \end{equation}
with the corresponding Hilbert product structure and we will consider the nonlinear mapping $H : Y \mapsto Z$, with
   \begin{equation} \label{defH}
H(y,v) := \left(y_t - \nabla\cdot\left[ a(\nabla y)\nabla y\right] - \chi_\omega v, y(0\,,\cdot)\right).
   \end{equation}
   In order to apply Theorem~\ref{Liusternik} in this setting, we must first establish some auxiliary results whose proofs are respectively given in Appendices~B, C and~D:
      
\begin{lemma} \label{WellDefiniteness}
   The mapping $H : Y \mapsto Z$ introduced in~\eqref{defH} is well-defined and continuous.
\end{lemma}

\begin{lemma} \label{ContinuousAndOnto}
   The linear mapping $\Lambda : Y \mapsto Z$, given by
   $$
\Lambda(y,v) = (y_t - a(0)\Delta y - \chi_\omega v, y(0\,,\cdot)) \quad 	\forall (y,v) \in Y,
   $$ 
is well-defined, continuous and onto.
   Furthermore, there exists a constant $M>0$ such that
   $$
\|\Lambda(y,v)\|_Z \geq M\|(y,v)\|_Y \quad \forall (y,v) \in Y.
   $$
\end{lemma}

\begin{lemma} \label{StrictDifferentiabilityAtTHeOrigin}
   The mapping $H$ is strictly differentiable at $(0,0) \in Y$ and~$H'(0,0) = \Lambda$.
\end{lemma}

\subsection{Proof of Theorem \ref{Maintheorem}}

   In view of Lemmas~\ref{WellDefiniteness} and~\ref{StrictDifferentiabilityAtTHeOrigin}, it is possible to apply Theorem~\ref{Liusternik} (for instance, with~$\epsilon = 1/(2M)$) and conclude that there exist~$\delta > 0$ and~$W : B_Z(0;\delta) \mapsto B_Y(0;2M\delta)$ such that
   \begin{equation} \label{RightInverseOfH}
H(W(f,y_0)) = (f,y_0) \quad \forall (f,y_0) \in B_Z((0,0);\delta).
   \end{equation}

   Now, if~$y_0 \in H^5(\Om) \cap H^1_0(\Om)$ with~$\|y_0\|_{H^5} \leq \delta$ and we take $(y,v) = W(0,y_0)$, we deduce from~(\ref{RightInverseOfH}) that
   \begin{equation} \label{ConclusionNCQuasi-linear}
(y_t - \nabla \cdot\left[ a(\nabla y) \nabla y \right] - \chi_\omega v, y|_{t=0}) = H(W(0,y_0)) = (0,y_0).   
   \end{equation}
   It is evident that (\ref{ConclusionNCQuasi-linear}) and the fact that~$(y,v) \in Y$ imply that $(y,v)$ solves~(\ref{sistema1fase}) and~(\ref{NC}), since the weights $\rho_k$ blow up to infinity exponentially as $t \to T$.

\subsection{Proof of Corollary \ref{corol-1.4}}

   Let $\gamma$ be given and let $\hat{\Om}$ be a bounded connected open set with a smooth boundary such that $\Om \subset \hat{\Om}$, $\partial\hat{\Om} \cap \partial\Om = \gamma$ and~$\hat{\omega} := \hat{\Om}\backslash \overline{\Om}$ is a non-empty open set.
   Let us set
   $$
\hat{Q} := ]0,T[ \times \hat{\Om}
\ \text{ and } \ \hat{\Sigma} := \left[0,T\right] \times \partial\hat{\Om}.
   $$
   There exists a linear continuous extension mapping $E : H^5(\Om) \cap H_0^1(\Om) \mapsto H^5(\hat{\Om}) \cap H_0^1(\hat{\Om})$, with
   $$
\| E y_0 \|_{H^5(\hat{\Om})} \leq \hat{C} \| y_0 \|_{H^5} \quad \forall y_0 \in H^5(\Om) \cap H_0^1(\Om).
   $$
   Then, we can apply Theorem~\ref{Maintheorem} to the system
   \begin{equation} \label{3.34p}
\begin{cases}
\hat{y}_t - \nabla\cdot\left[ a(\nabla \hat{y}) \nabla \hat{y}\right] = \chi_{\hat{\omega}} \hat{v} \ &\text{ in } \ \hat{Q}, \\
\hat{y} = 0 \ &\text{ on } \ \hat{\Sigma}, \\
\hat{y}|_{t=0} = \hat{y}_0 \ &\text{ in }\ \hat{\Om},
\end{cases}
   \end{equation}
whence we deduce that there exists $\eta > 0$ such that, if $\| \hat{y}_0 \|_{H^5(\hat{\Om})} \leq \eta$, \eqref{3.34p} admits a control that drives the state exactly to zero at~$t = T$.

   Consequently, if $y_0$ satisfies $\| y_0 \|_{H^5} \leq \delta := \eta/\hat{C}$ and~$y$ and~$h$ are respectively given by the restriction to~$Q$ and the trace on~$]0,T[ \times \gamma$ of the corresponding $\hat{y}$, we find that~$(y,h)$ is a state-control pair that solves the considered boundary control problem.

\section{The convergence of \textbf{ALG~1}, approximations and experiments}\label{Sec4}

   As we already explained in Section~\ref{intro}, arguing  as in~\cite{CCLM, CHCV}, we can introduce an elementary quasi-Newton algorithm for the computation of a solution to the null control problem.
   To this purpose, we previously have to specify the inverse to the right of~$H'(0,0)$.
   We can do this with the help of the Fursikov-Imanuvilov method~\cite{Fursikov-Imanuvilov}.

   The argument is the following.
   For any $(f,y_0) \in Z$, we can obtain a solution to~\eqref{Sistema1FaseLinear} in~$Y$ by solving the following extremal problem:
    \begin{equation}
  \label{eq.43'}
  \left\lbrace
  \begin{array}{ll}
  \displaystyle \mbox{Minimize } \jjntQ \rho^2 |y|^2 + \jjntomT  \rho_{3}^2 |v|^2 \\
  \mbox{Subject to } v \in L^2(]0,T[ \times \omega), \ (y,v) \mbox{ satisfies } \eqref{sistema1fase}.
  \end{array}
  \right.
    \end{equation}

   Recall that the definition of the weights $\rho_k$ is given in~\eqref{Sistema1FaseLinear}.

   It is well-known that \eqref{eq.43'} possesses exactly one solution, given by
   \begin{equation}
\label{eq.43''}
 y=\rho^{-2} L^*p, \quad \widehat{v} = -\rho^{-2}_3 p, \quad v= \left. \chi \widehat{v}\right|_{]0,T[ \times \omega},
   \end{equation}
where $p$ is the unique solution to the variational equality
   \begin{equation}
\label{eq.43'''}
\left\lbrace
   \begin{array}{ll}
   \displaystyle \pi(p,w)= \jnt_\Om y_0(x) w(0,x) \,dx + \jjntQ f w\\
   \forall w \in P , \ \ p \in P
   \end{array}
   \right.
   \end{equation}
   (recall the notation introduced in the proof of~Theorem~\ref{Maintheorem}).
   The existence and uniqueness of a solution to~\eqref{eq.43'''} is an immediate consequence of Lax-Milgram Theorem; see the details in~\cite{Fursikov-Imanuvilov}.

   Accordingly, we set $H'(0,0)^{-1} (f, y_0) = (y,v)$, with $y$ and~$v$ respectively given by~\eqref{eq.43''}--\eqref{eq.43'''}.

   Note that, in \textbf{ALG~1}, for each $n$ the task reduces to solve a problem of the kind \eqref{eq.43'''}, with
   \begin{equation}\label{def_f}
f = - \nabla \cdot ((a(0)  - a(\nabla y^n)) \nabla y^n)).
   \end{equation}
   The local convergence of \textbf{ALG~1} is guaranteed by the following result:

\begin{theorem}\label{th.1.3}
   Let $y_0 \in H^5(\Om) \cap H^1_0(\Om)$ be given, with $\Vert y_0 \Vert_{H^5} \leq \varepsilon$, where $\varepsilon$ is the small quantity furnished by Theorem~\ref{Maintheorem} ensuring the null controllability of~\eqref{sistema1fase}.
   There exists $\kappa > 0$ such that, if $(y^0, v^0)\in Y$ and
   \[
\Vert (y^0,v^0) - (y,v) \Vert_Y \leq \kappa,
   \]
then we get linear rate convergence to~$(y,v)$ in~$Y$ for the sequence $\left\{ (y^n, v^n) \right\}_n$.
\end{theorem}

   The proof of this result is very similar to the proof  of~Theorem~1.2 in~\cite{fernandez2021theoretical}.
   For brevity, it is omitted.

   In order to solve numerically the problems \eqref{eq.43'''}, a first idea is to construct adequate finite dimensional spaces $P_h \subset P$ and search for functions $p_h \in P_h$ satisfying the equalities in~\eqref{eq.43'''} for all $w \in P_h$.

   This is relatively easy when~$d = 1$.
   However, if~$d \geq 2$, this is not so simple.
   The reason is that the functions in $P$ must satisfy $L^* p \in L^2_{loc}(Q)$ and, consequently, an approximation based on a standard mesh of $Q$ requires spaces $P_h$ of functions that must be globally $C^0$ in all the variables and globally $C^1$ in space.
    Accordingly, for $d \geq 2$, it seems convenient to use a mixed formulation, as  in~\cite{FCMS, CHCV}.

   For completeness, we have performed numerical experiments for~$d=1$ with approximations of both kinds (based on subspaces~$P_h$ and also based on mixed formulations).
   For $d=2$, our tests only concern mixed finite elements.

   The details are given in~\cite{FCMunch, FCMS}; see also~\cite{FCMunchNewton}.

\subsection{A direct approximation of~\eqref{eq.43'''} for $d=1$: $\mathbb{Q}_3-\mathbb{Q}_1$ finite elements}

   As already mentioned, it is not difficult to introduce efficient finite dimensional subspaces of~$P$ in this case.
   We take $\Om = ]0,L[$ and consequently have $Q = ]0,T[ \times ]0,L[$ and~$\Sigma = ]0,T[ \times \{0,L\}$.
   Let $N_x$ and~$N_t$ be two (large) positive integers and set $\Delta x = L/N_x$, $\Delta t = T/N_t$ and~$h = (\Delta x,\Delta t)$.
   Then, we replace in~\eqref{eq.43'''} the space $P$ by the finite dimensional subspace $P_h$ formed by the functions $p_h \in C^0(\overline{Q})$ such that $p_{h,x} \in C^0(\overline{Q})$, the restriction to each rectangle is polynomial of degree $\leq 3$ in~$x$ and degree $\leq 1$ in~$t$ and vanish on~$\Sigma$.

   This approach was introduced in~\cite{FCMunch}.
   A complete analysis of the convergence of the resulting approximations was given there.

\begin{table}[htbp]
\centering
\caption{Test 1.1 -- The control and state relative errors with respect to the ``exact'' solution corresponding to the solutions computed for various $h$.}
\label{Table_1.1-1}
\begin{tabular}{cccc}
\toprule
$h = (\Delta x,\Delta t)$  & Number of Points  & Control rel error & State rel error \\
\midrule
$ (0.1,0.05) $            &  $100$   & $0.4170$ & $0.2252$\\
$ (0.05,0.025) $        &  $400$   & $0.0265$ & $0.0214$\\
$ (0.0333,0.0167) $  &  $900$   & $0.0037$ & $0.0033$\\
$ (0.025,0.0125) $    &  $1600$ & $0.0012$ & $0.0010$\\
$ (0.02,0.01) $          &  $2500$ & $0.0005$ & $0.0004$\\
$ (0.0167,0.0083) $  &  $3600$ & --- & --- \\
\bottomrule
\end{tabular}
\end{table}

\subsection{A mixed finite dimensional approximation of~\eqref{eq.43'''} for $d=1$  and~$d=2$}

    In order to get a formulation where the finite dimensional spaces are more tractable, we can follow a different process:

\begin{enumerate}

\item First, we introduce new variables $(z,m)$ and we rewrite~\eqref{eq.43'''} accordingly.
   Indeed, since the weights $\rho$ and $\rho_3$ blow up as~$t \to T$, in practice it is convenient to work with the new variables~$z= \rho^{-1} L^*p$ and~$m= \rho_{3}^{-1} p$.
   We make use of the equality $z= \rho^{-1} L^*(\rho_3 m)$ to introduce a multiplier and get a mixed formulation.

\item Then, we integrate by parts and rewrite the formulation in such a way that only first order partial derivatives remain.

\item Finally, we replace the resulting mixed formulation by a finite dimensional approximation using standard Lagrange ($C^0$ in time and space) finite elements.

\end{enumerate}

   This mixed formulation was also introduced in~\cite{FCMunch} and then reconsidered and adapted to several situations on other papers, see for instance~\cite{FCMS}.

   In the following sections, we will provide the results of some experiments.

\subsection{Test 1.1}

   We applied the quasi-Newton method to the solution to the null controllability problem for~\eqref{sistema1fase} with the following data:
\begin{itemize}
\item $d=1$, $\Om = ]0,1[$, $\omega = ]0.3, 0.7[$, $T = 0.5$.
\item $y_0(x) \equiv 0.9 \sin(\pi x)$.
\item $a(s) \equiv 1.5+ 0.5 \sin(s)$.
\end{itemize}

   At each step of \textbf{ALG~1}, we wrote the problem in the form \eqref{eq.43'''}.
   This first test corresponds to experiments performed by direct approximation (Sect.~4.1) with different values of~$h = (\Delta x,\Delta t)$.
   The computations and visualizations have been carried out with the help of the MatLab Library~\cite{MATLAB:2010}.

   In this and the other tests, we  fix the stopping criterion
   \begin{equation}\label{stopping}
\displaystyle \frac{\Vert y^{n+1} - y^{n}\Vert_{L^2(Q)}}{\Vert y^{n+1} \Vert_{L^2(Q)}} \leq \kappa ,
   \end{equation}
where $y^n$ and~$y^{n+1}$ are the computed states respectively at the $n$-th and~$(n+1)$-th steps and~$\kappa = 10^{-5}$.

   The evolution of the  relative errors of the computed solutions for various $h$ is given in Table~\ref{Table_1.1-1}.
   There, we have assigned the role of ``exact'' to the solution computed for $\Delta x = 0.0167$ and~$\Delta t = 0.083$, which corresponds to a mesh with $3600$ points.

\begin{figure}[htbp]
\centering
\includegraphics[height=52mm, width = 62mm]{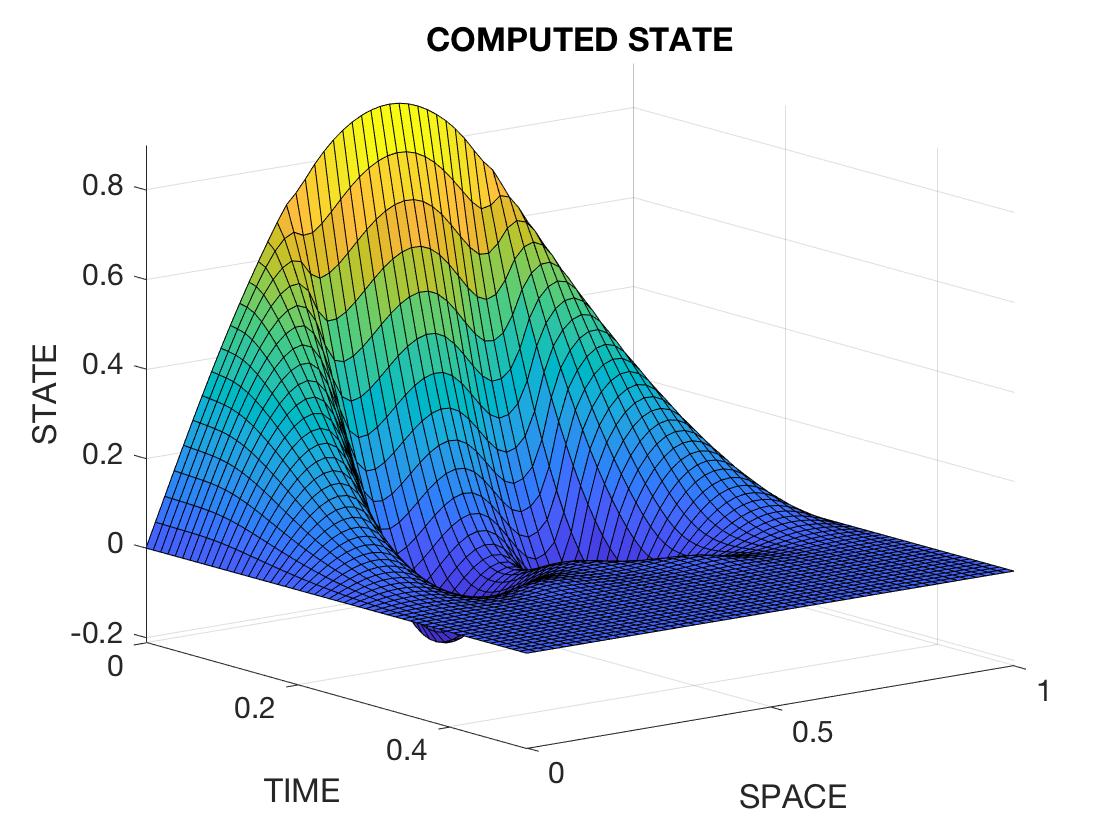}\quad
\includegraphics[height=52mm, width = 62mm]{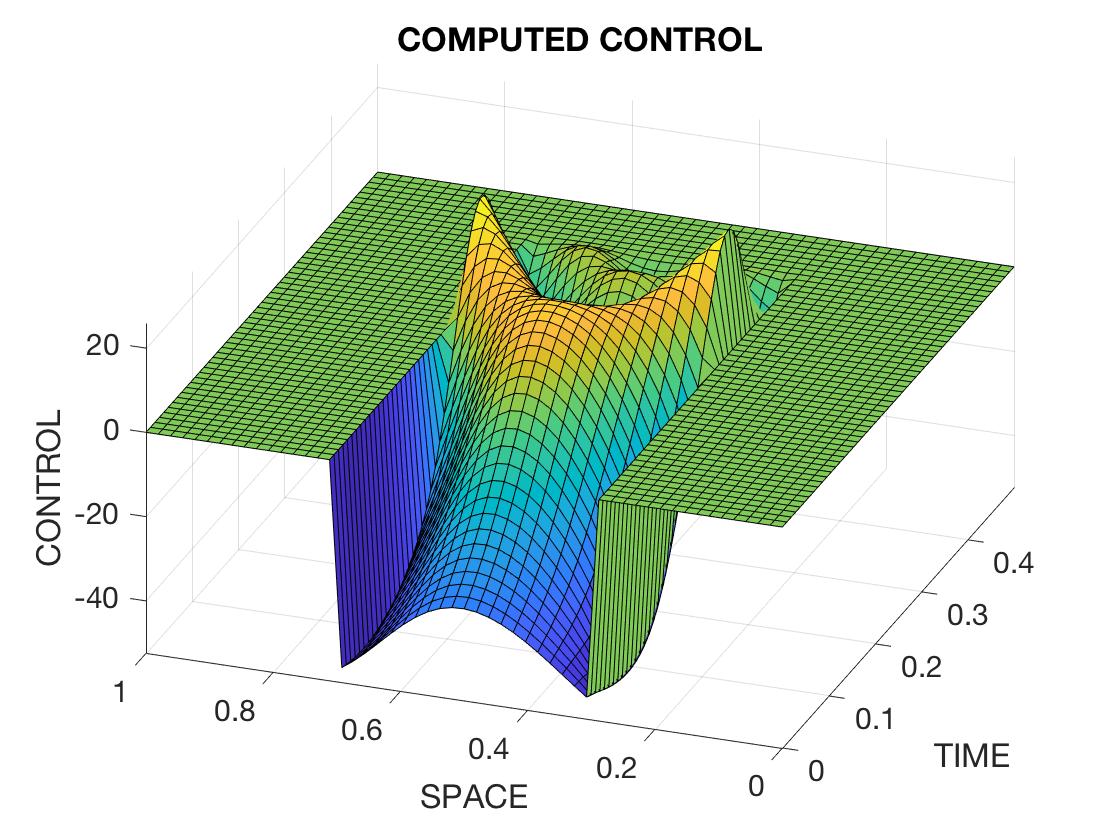}
\caption{Test 1.1 -- The computed state and control for~$\Delta x = 0.02$ and~$\Delta t = 0.01$.}
\label{Fig_1.1-1}
\end{figure}

\begin{figure}[htbp]
\centering
\includegraphics[height=48mm, width = 60mm]{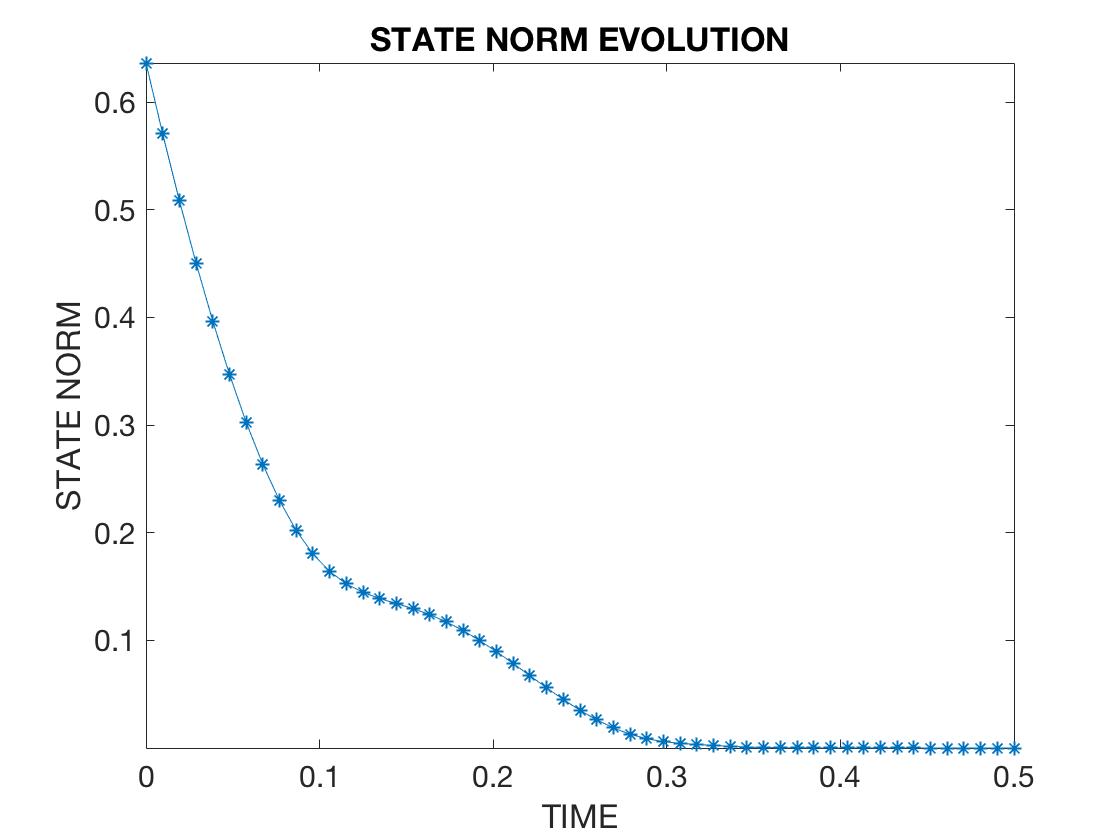}\quad
\includegraphics[height=48mm, width = 60mm]{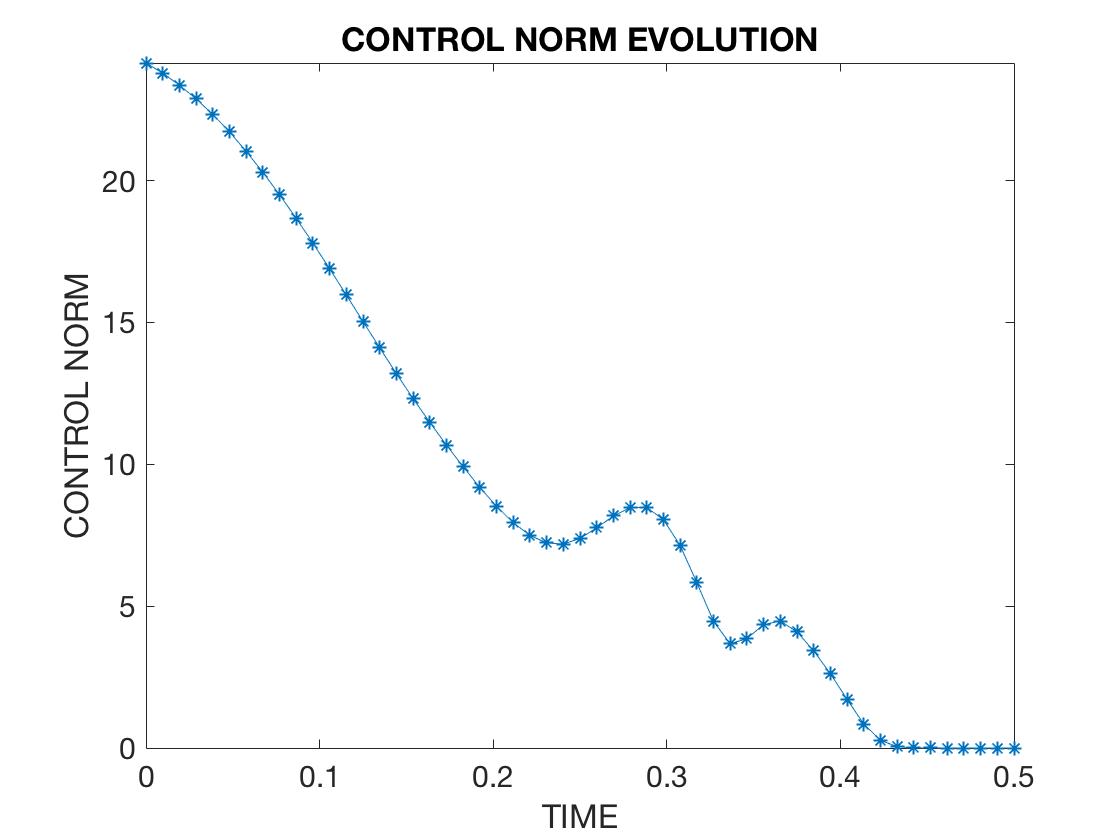}
\caption{Test 1.1 -- Evolution in time of the $L^2$ norms of the state and the control for~$\Delta x = 0.02$ and~$\Delta t = 0.01$.}
\label{Fig_1.1-2}
\end{figure}

   The computed state and control for~$\Delta x = 0.02$ and~$\Delta t = 0.01$ are depicted in~Figure~\ref{Fig_1.1-1}.
   The evolution in time of the associated $L^2$ norms is given in~Figure~\ref{Fig_1.1-2}.

\subsection{Test 1.2}

   We take the same data except for the initial state and the diffusion coefficient.
   This time,

\begin{itemize}
\item $y_0(x) \equiv 0.15 \sin(\pi x)$.
\item $a(s) \equiv \exp(-2 \exp(-0.3s))$.
\end{itemize}

\begin{figure}[htbp]
\centering
\includegraphics[height=52mm, width = 62mm]{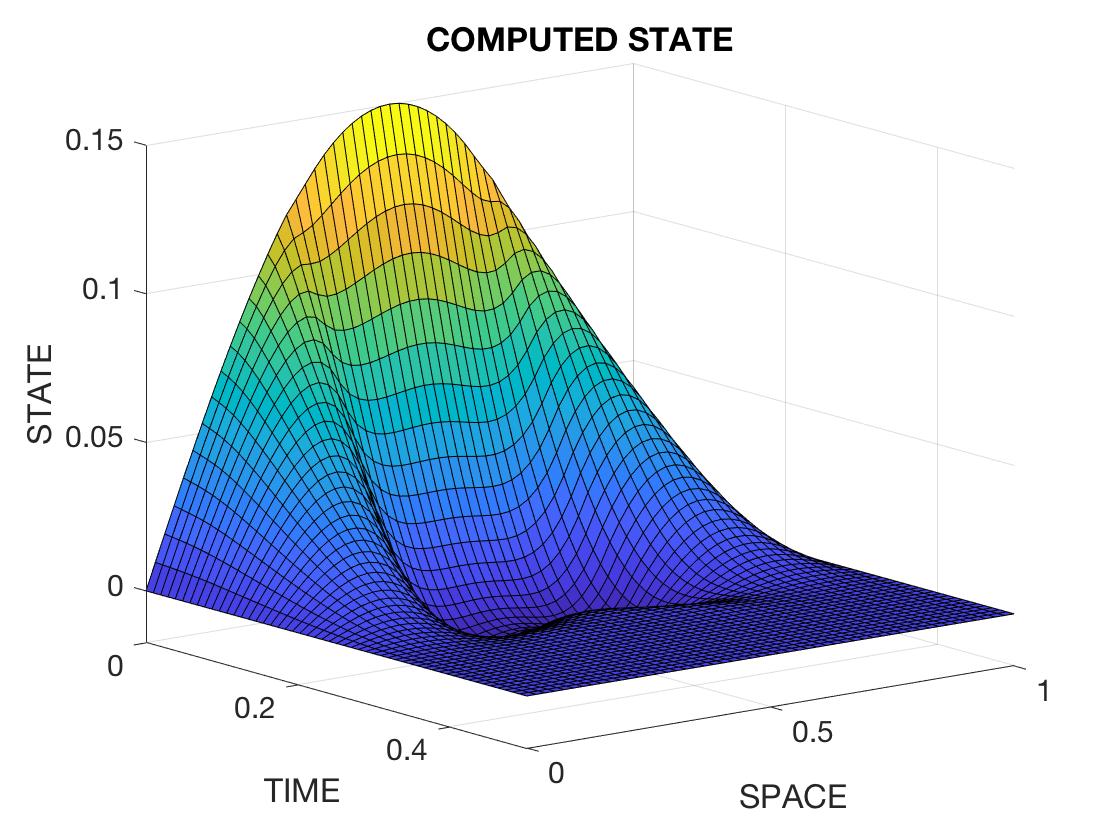}\quad
\includegraphics[height=52mm, width = 62mm]{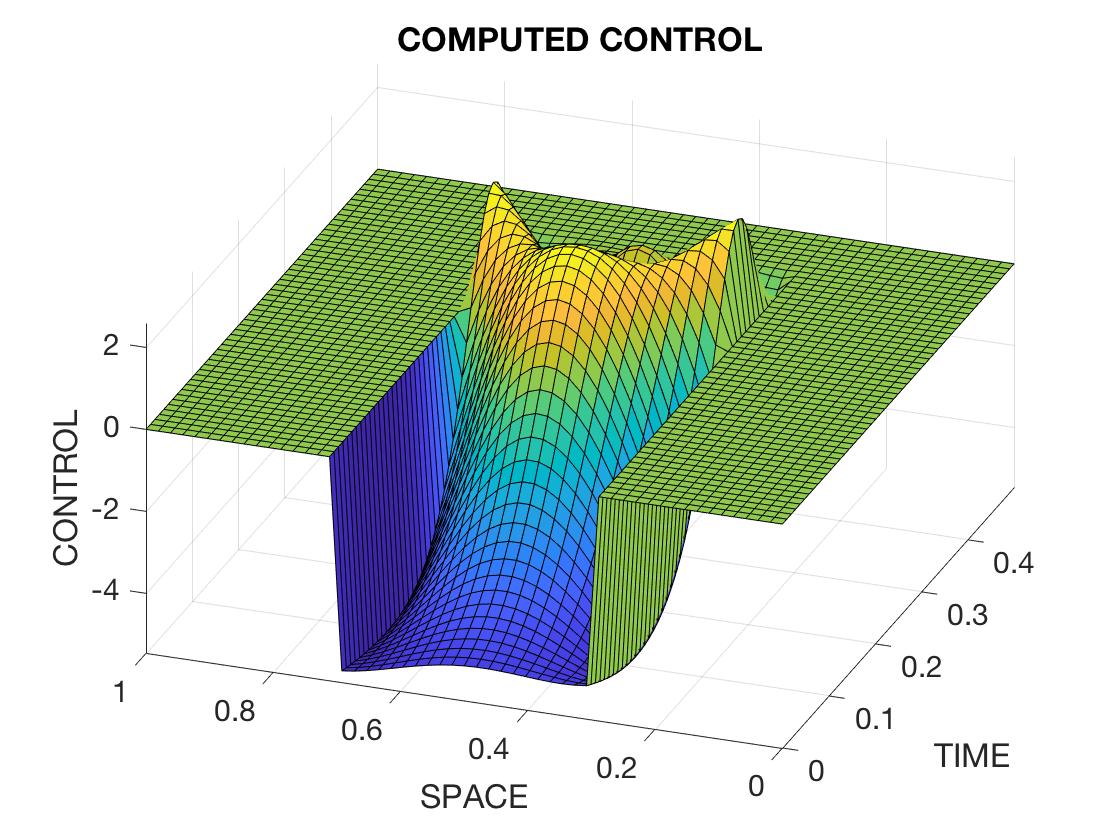}
\caption{Test 1.2 -- The computed state and control for~$\Delta x = 0.02$ and~$\Delta t = 0.01$.}
\label{Fig_1.2-1}
\end{figure}

\begin{figure}[htbp]
\centering
\includegraphics[height=48mm, width = 60mm]{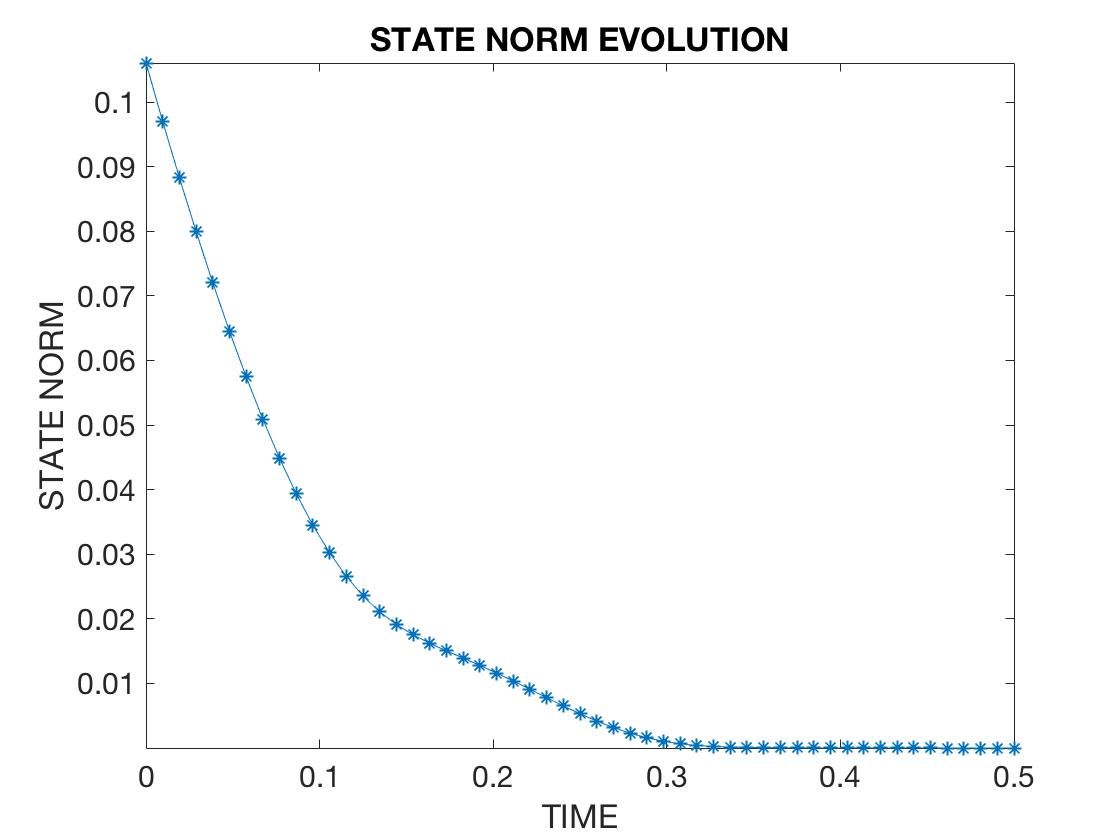}\quad
\includegraphics[height=48mm, width = 60mm]{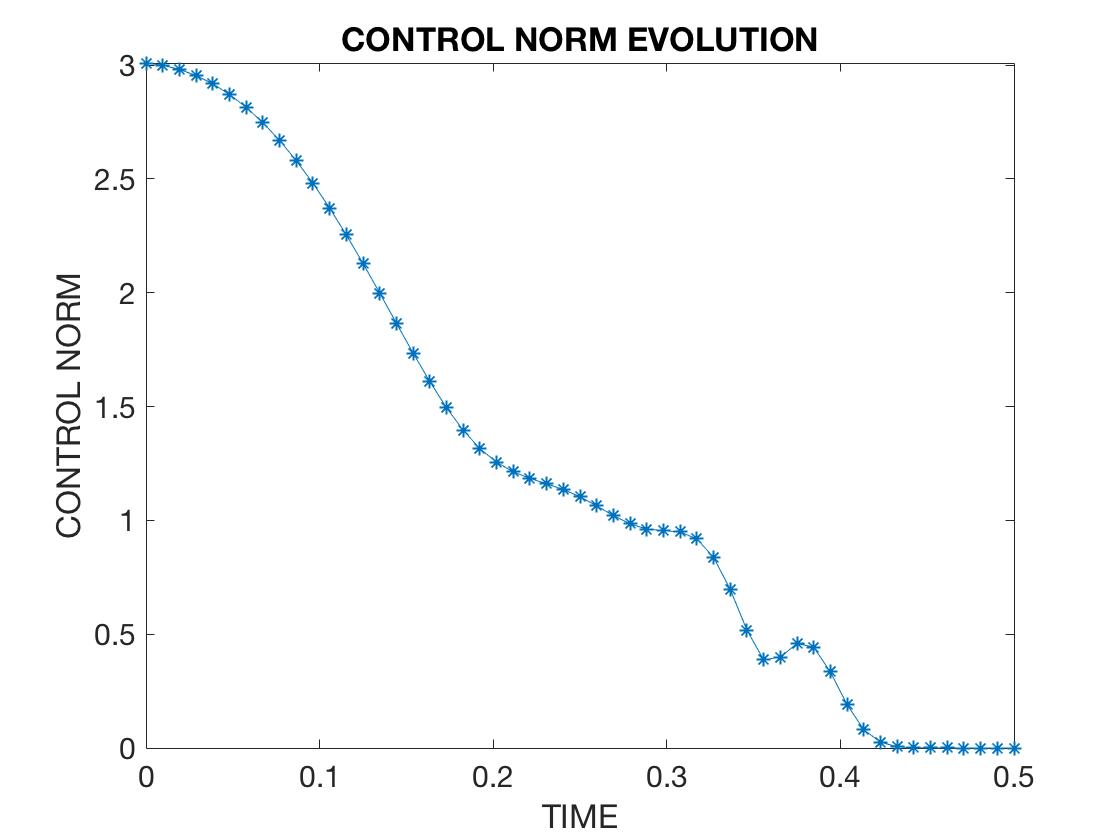}
\caption{Test 1.2 -- Evolution in time of the $L^2$ norms of the state and the control for~$\Delta x = 0.02$ and~$\Delta t = 0.01$.}
\label{Fig_1.2-2}
\end{figure}

   Again, we approximate~\eqref{eq.43'''} as in Section~4.1 with the same values of $h = (\Delta x,\Delta t)$.
   The computed state and control and the evolution in time of the corresponding norms is now given in~Figures~\ref{Fig_1.2-1}--\ref{Fig_1.2-2}.

   A comparison of these results with those in Test~1.1 shows that, although the initial state is larger, the norm of the computed control is smaller.
   This seems to be the effect of a monotone diffusion coefficient.
   However, this should be investigated more in detail.

\subsection{Test 1.3}

   In this test, we take again the data in Test~1.1.
   At each step, we compute a numerical approximation of the corresponding problem~\eqref{eq.43'''} using piecewise linear and quadratic mixed finite elements (Sect.~4.2).

   We have performed the computations with the {\tt FreeFem++} package, with a mesh adaptation technique.
   This means that, at each quasi-Newton iterate, the mesh is not regular, but adapted to the behavior of the last computed solution.
   More precisely, the points are distributed in the space-time domain in accordance with a metric that corresponds to the values of the computed (approximate) Hessian of~$p^n$;
   related explanations are given in~\cite{Hecht-adapt}.
   This technique was already used, for instance, in reference~\cite{fernandez2000null} for the solution of another null controllability problem; see Figure~2 there.

   For a detailed description of {\tt FreeFem++,} see~\cite{FreeFem-1, FreeFem-2}.

\begin{figure}[htbp]
\centering
\includegraphics[width= 72mm, height= 72mm]{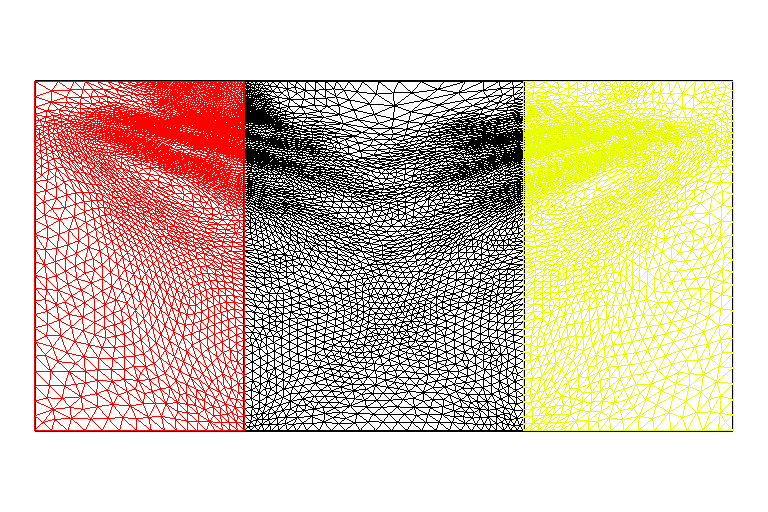}
\caption{Test 1.3 -- The space-time domain is $]0,1[\times]0,0.5[$; the control domain is the central band $]0.3,0.7[\times]0,0.5[$. The final mesh after adaptation is shown. Number of vertices: 6189. Number of triangles: 12157.}
\label{Fig_1.3-2}
\end{figure}

\begin{figure}[htbp]
\centering
\includegraphics[height=52mm, width = 62mm]{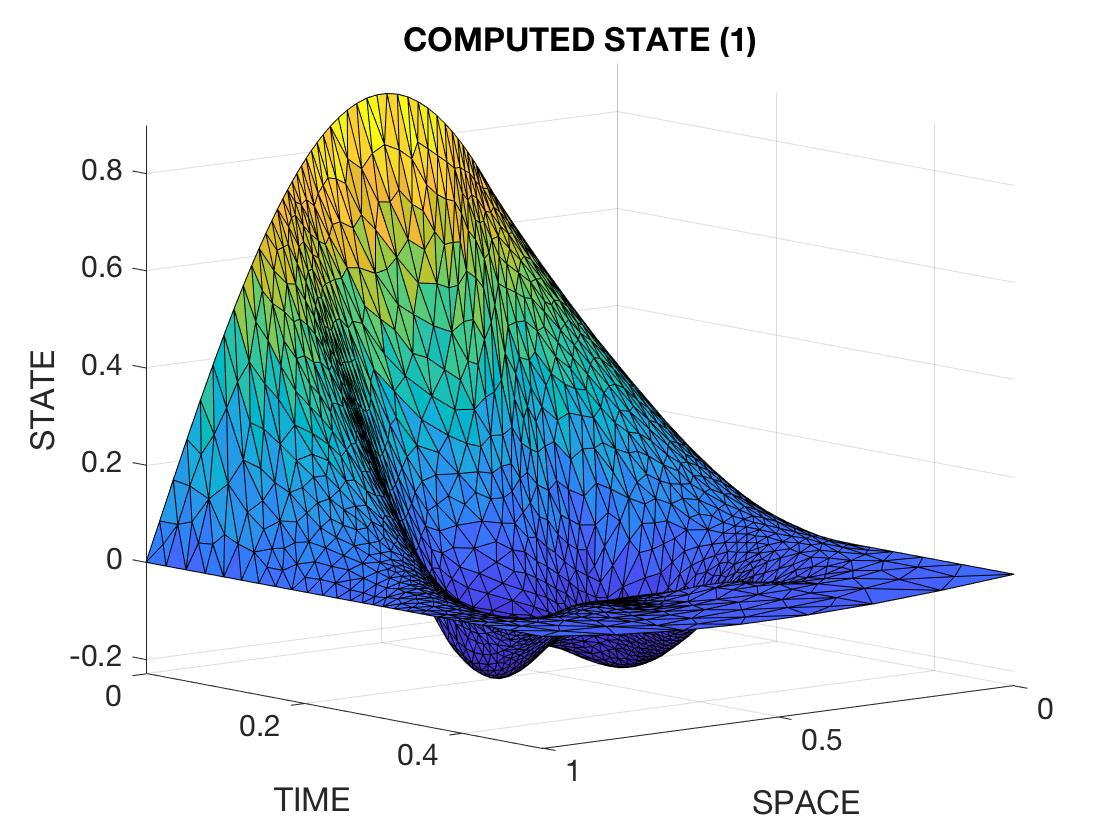}\quad
\includegraphics[height=52mm, width = 62mm]{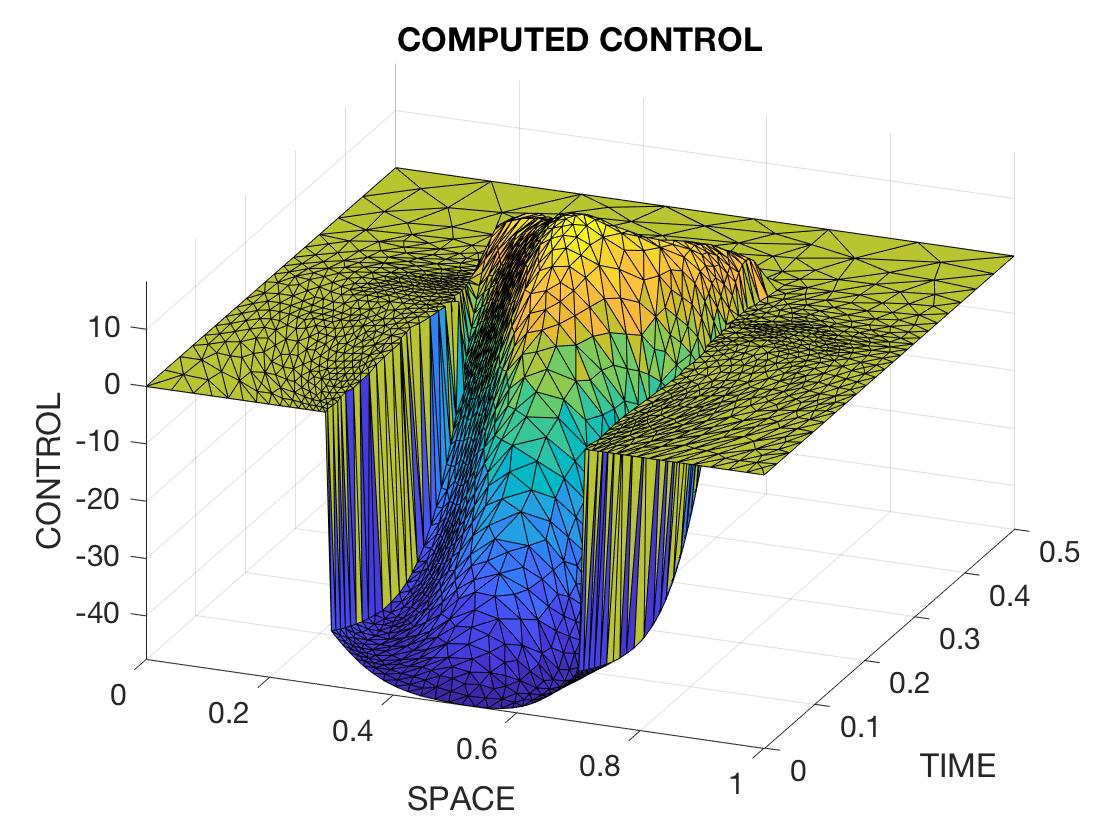}
\caption{Test 1.3 -- The computed state and control. The state (on the left) and the control (on the right) are found by solving numerically~\eqref{eq.43'''} with mixed finite element techniques.}
\label{Fig_1.3-3}
\end{figure}

   The stopping criterion is~\eqref{stopping}, where the notation is self-explanatory.
   We have started from a regular mesh with $970$ vertices and $1844$ triangles.
   The final mesh is displayed in~Figure~\ref{Fig_1.3-2}.
   The computed state and control are displayed in~Figure~\ref{Fig_1.3-3}.

   They have been compared to the most accurate results in~Test~1.1 in~Table~\ref{Table_1.3-1}, with excellent agreement.

\begin{figure}[htbp]
\centering
\includegraphics[height=52mm, width = 62mm]{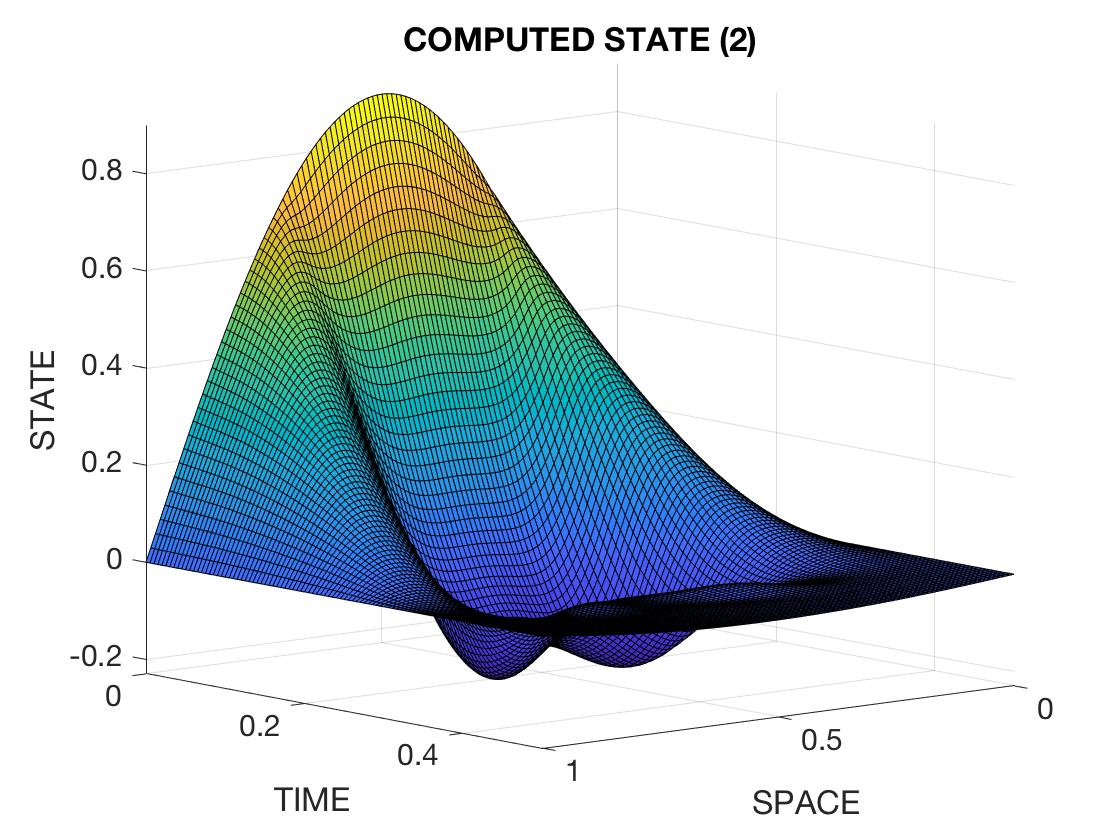}\quad
\includegraphics[height=52mm, width = 62mm]{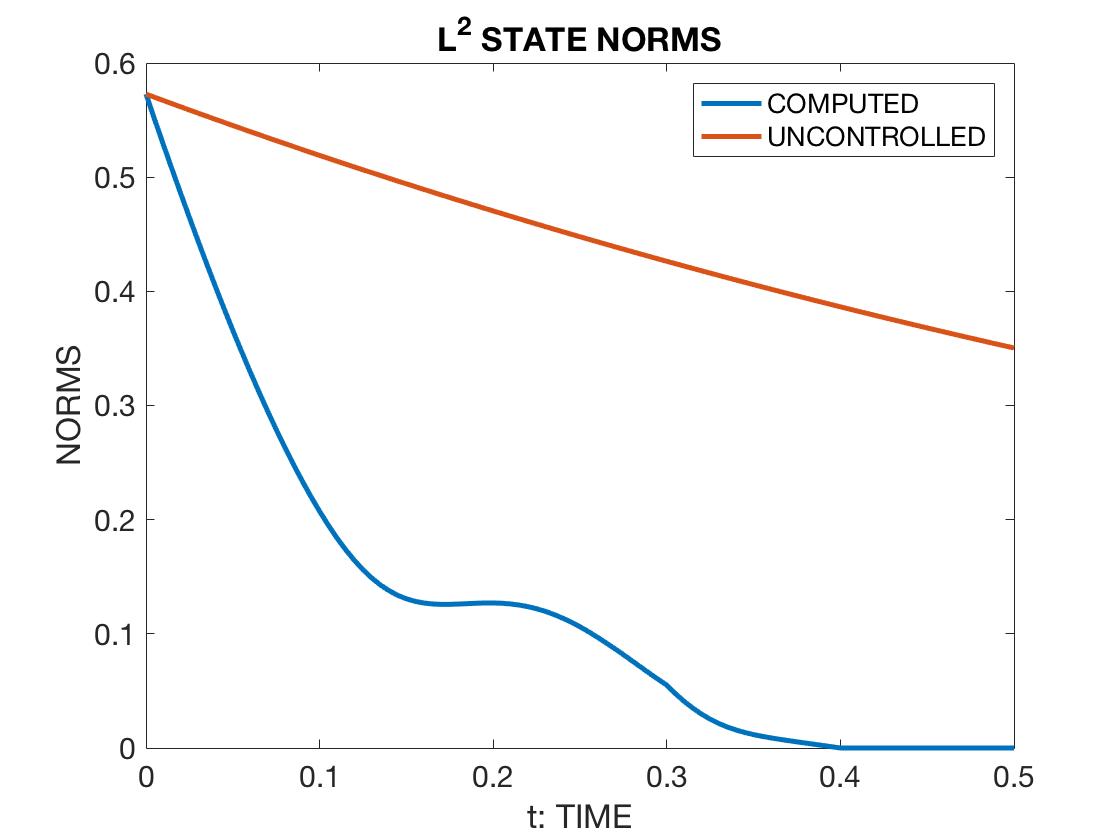}
\caption{Test 1.3 -- The computed state and the norm evolution. Here, the state (on the left) is found by solving numerically~\eqref{sistema1fase}. On the right, a comparison of the spatial $L^2$-norms of the computed state and the solution to~\eqref{sistema1fase} for $v \equiv 0$.}
\label{Fig_1.3-3Bis}
\end{figure}

   Also, in order to check that the computed control does the work, we have solved numerically the nonlinear problem~\eqref{sistema1fase} with this control in the right hand side.
  To this end, we have used the MatLab tool {\tt pdepe,} that uses {\it second order} finite difference approximation in space and a variable-step, {\it fifth-order} time discretization scheme based on numerical differentiation, see~\cite{Skeel}.
  The state found this way is depicted in~Figure~\ref{Fig_1.3-3Bis} and is in practice identical to the solution to the test.
  This shows that the method described in Section~4.2 certainly furnishes a null control.

\begin{table}[htbp]
\centering
\caption{Tests 1.1 and 1.3 -- Comparisons of the best computed state and control in Test 1.1 and the computed state and control in Test 1.3.}
\label{Table_1.3-1}
\begin{tabular}{cccc}
\toprule
Control $L^2$ norm (Test 1.1) & Control $L^2$ norm (Test 1.3) & Norm of difference & Rel difference \\
\midrule
$42.0693$  & $42.0488$ & $0.0350$ & $0.0004$ \\
\midrule
State $L^2$ norm (Test 1.1) & State $L^2$ norm (Test 1.3) & Norm of difference & Rel difference \\
\midrule
$0.6726$  &  $0.6717$ & $0.0008$ & $0.0012$ \\
\bottomrule
\end{tabular}
\end{table}

   In the framework of this test, we have studied how large the initial data can be to get convergence for \textbf{ALG~1}.
   Specifically, we have taken
 \begin{itemize}
\item $y_0(x) \equiv y_{00} \sin(\pi x)$, with several $y_{00}$.
\end{itemize}
   The amount of iterates needed to ensure~\eqref{stopping} can be found in~Table~\ref{Table_1.3-2}.
   Convergence fails for $y_{00} \geq 2.06$;
   it is reasonable to expect that this corresponds to initial data too far from zero to support \textbf{ALG~1} and, maybe, a continuation algorithm is in order.\footnote{This was suggested by one of the referees.}

   We have also compared \textbf{ALG~1} and the usual Newton method, for which the iterates are given by
   $$
(y^{n+1}, v^{n+1}) =(y^n, v^n) - H'(y^n, v^n)^{-1} (H(y^n, v^n) - (0,y_0)).
   $$
   In practice, this means that, at each step, we must solve a problem like~\eqref{eq.43'''} with a different $\pi(\cdot\,,\cdot)$ and $f$ is as in~\eqref{def_f}.
   The results of this comparison are given in~Table~\ref{Table_1.3-3}.
   It is readily seen that \textbf{ALG~1} is clearly preferable both in terms of CPU execution time and similar in terms of convergence rate.

\begin{table}[htbp]
\centering
\caption{Tests 1.3 -- Iterates for convergence vs.\ initial data size. The first and third (resp.\ second and fourth) columns indicate the numbers of iterates needed to satisfy the stopping test for $\kappa = 10^{-5}$ (resp.\ $\kappa = 10^{-9}$). Convergence fails for $y_{00} \geq 2.06$.}
\label{Table_1.3-2}
\begin{tabular}{cccccc}
\toprule
$y_{00}$ & Iter.\ for $\kappa=10^{-5}$ & Iter.\ for $\kappa=10^{-9}$ & $y_{00}$ & Iter.\ for $\kappa=10^{-5}$
& Iter.\ for $\kappa=10^{-9}$ \\
\midrule
$0.05$ & $4$ & $6$ & $1.05$ & $25$ & $43$ \\
$0.15$ & $6$ & $11$ & $1.30$ & $29$ & $50$ \\
$0.30$ & $9$ & $15$ & $1.45$ & $29$ & $50$ \\
$0.45$ & $10$ & $15$ & $1.60$ & $28$ & $50$ \\
$0.60$ & $10$ & $16$ & $1.75$ & $29$ & $50$ \\
$0.75$ & $11$ & $16$ & $1.90$ & $39$ & $50$ \\
$0.90$ & $21$ & $37$ & $2.05$ & $206$ & $372$ \\
\bottomrule
\end{tabular}
\end{table}

\begin{table}[htbp]
\centering
\caption{Tests 1.3 -- Comparison of \textbf{ALG~1} (QN iterates, rate, etc.) and Newton algorithm (Nw iterates, Nw rate, etc.). It is observed that \textbf{ALG~1} and Newton algorithm are similar in convergence rate but \textbf{ALG~1} is much less expensive. An explanation is that the domain of convergence of \textbf{ALG~1} is larger and increasing initial data slow down Newton iterates before. The CPU time is given in seconds.}
\label{Table_1.3-3}
\begin{tabular}{ccccccc}
\toprule
$y_{00}$ & QN iterates & QN rate & QN cpu time & Nw iterates & Nw rate & Nw cpu time \\
\midrule
 $0.050000$ & $7$ & $1.13541$ & $0.052182$ & $5$ & $1.19228$ & $0.665083$ \\
 $0.057881$ & $7$ & $1.15098$ & $0.052295$ & $6$ & $1.15685$ & $0.665083$ \\
 $0.067005$ & $8$ & $1.12136$ & $0.051693$ & $6$ & $1.15660$ & $0.670934$ \\
 $0.077566$ & $8$ & $1.12116$ & $0.054900$ & $7$ & $1.15789$ & $0.678381$ \\
 $0.089793$ & $9$ & $1.11714$ & $0.052194$ & $7$ & $1.15737$ & $0.673889$ \\
 $0.103946$ & $9$ & $1.11765$ & $0.052057$ & $8$ & $1.11852$ & $0.664614$ \\
 $0.120331$ & $10$ & $1.09586$ & $0.053173$ & $9$ & $1.11677$ & $0.665924$ \\
 $0.139298$ & $10$ & $1.09565$ & $0.055559$ & $10$ & $1.09720$ & $0.668232$ \\
 $0.161255$ & $11$ & $1.09705$ & $0.052428$ & $11$ & $1.09307$ & $0.673780$ \\
 $0.186673$ & $12$ & $1.07911$ & $0.053231$ & $12$ & $1.08371$ & $0.665457$ \\
 $0.216097$ & $13$ & $1.08275$ & $0.053231$ & $14$ & $1.07392$ & $0.672346$ \\
 $0.250159$ & $14$ & $1.06756$ & $0.051752$ & $15$ & $1.07029$ & $0.678380$ \\
 $0.319274$ & $15$ & $1.07319$ & $0.051525$ & $16$ & $1.06058$ & $0.665457$ \\
 $0.369599$ & $15$ & $1.07247$ & $0.053866$ & $15$ & $1.06939$ & $0.671037$ \\
 $0.427858$ & $15$ & $1.07014$ & $0.049892$ & $14$ & $1.06443$ & $0.700798$ \\
 $0.495299$ & $15$ & $1.06997$ & $0.047642$ & $18$ & $1.04828$ & $0.689229$ \\
 $0.573370$ & $15$ & $1.06952$ & $0.052191$ & $77$ & $1.00933$ & $0.707382$ \\
 $0.663747$ & $15$ & $1.05980$ & $0.051467$ & $31$ & $1.02951$ & $0.665708$ \\
 $0.768371$ & $15$ & $1.08188$ & $0.047299$ & $77$ & $1.00933$ & $0.681966$ \\
 $0.806789$ & $17$ & $1.07032$ & $0.054590$ & $189$ & $1.00875$ & $0.686459$ \\
 $0.847129$ & $21$ & $1.05159$ & $0.054594$ & $196$ & $1.00589$ & $0.673658$ \\
 $0.889485$ & $27$ & $1.03753$ & $0.054181$ & -- & -- & -- \\
\bottomrule
\end{tabular}
\end{table}

\subsection{Test 2.1}

   The data in this Test are the following:

\begin{itemize}
\item $d=2$, $\Om = ]0,1[ \times ]0,1[$, $\omega = ]0.2, 0.8[ \times ]0.2, 0.8[$, $T = 0.5$.
\item $y_0(x_1,x_2) \equiv 0.4 \sin(\pi x_1) \sin(2\pi x_2)$.
\item $a(s) \equiv \exp(-2 \exp(-(s_1 + s_2))$.
\end{itemize}

\begin{figure}[htbp]
\centering
\includegraphics[width= 82mm, height= 70mm]{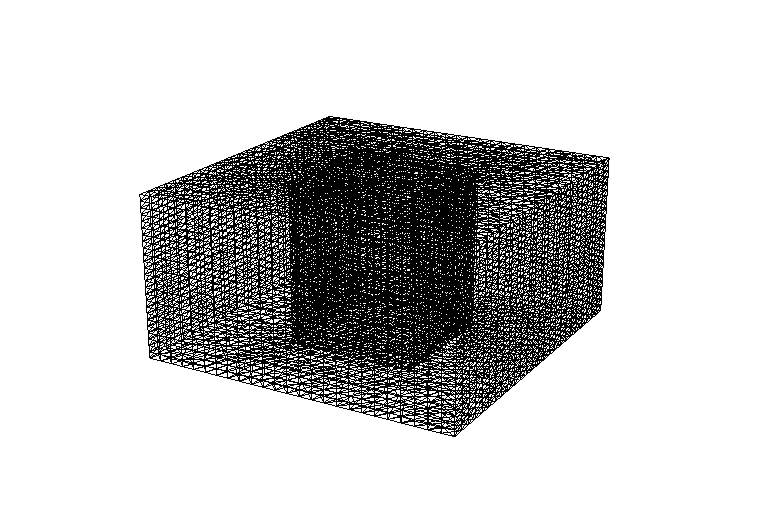}
\caption{Test 2.1 -- The finest mesh. Number of vertices: 20482, number of tetrahedra = 111,888.}
\label{Fig_2.1-1}
\end{figure}

\begin{figure}[htbp]
\centering
\includegraphics[height=48mm, width = 62mm]{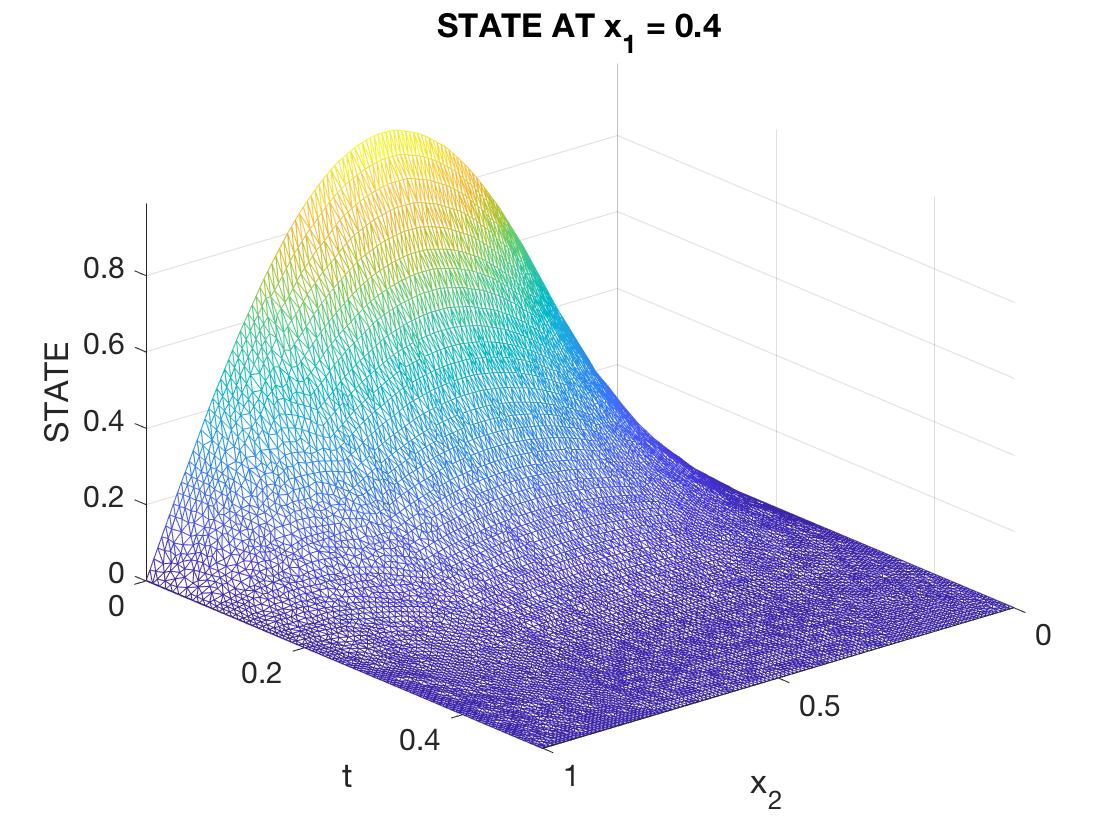}\quad
\includegraphics[height=48mm, width = 62mm]{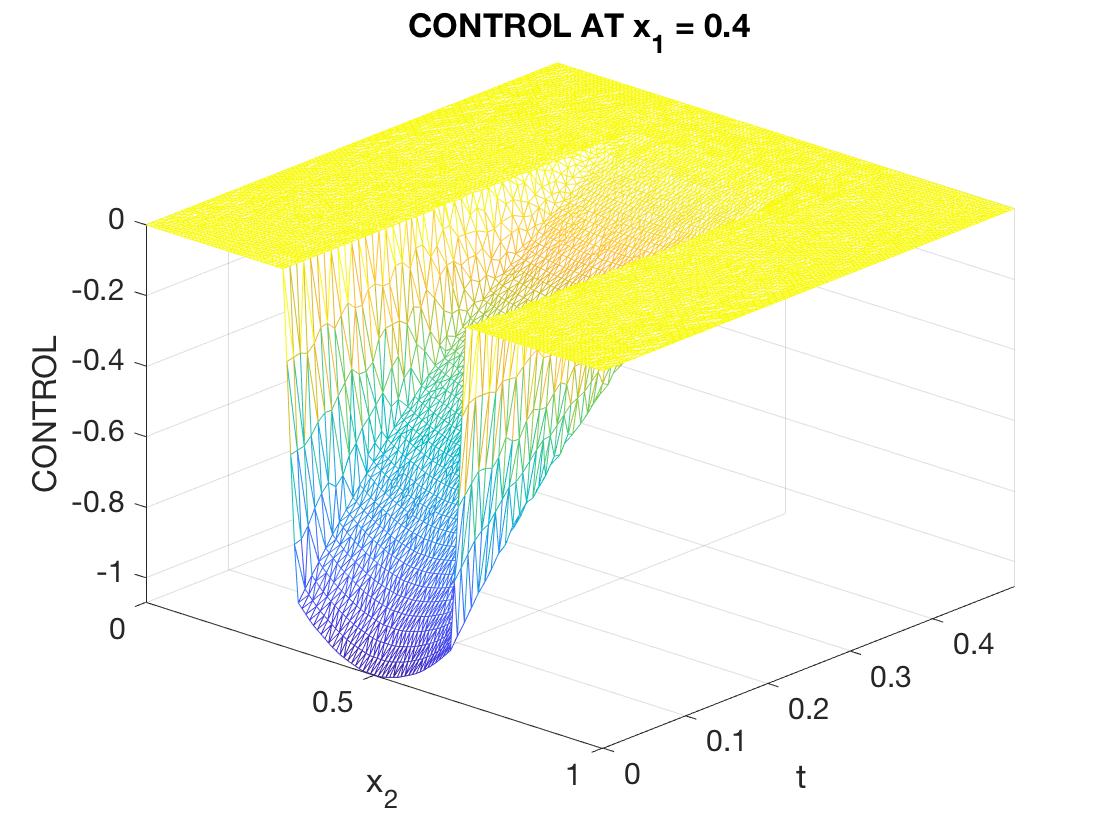}
\caption{Test 2.1 -- Cuts of the computed state and control at $x_1 = 0.4$.}
\label{Fig_2.1-3}
\end{figure}

\begin{figure}[htbp]
\centering
\includegraphics[height=48mm, width = 62mm]{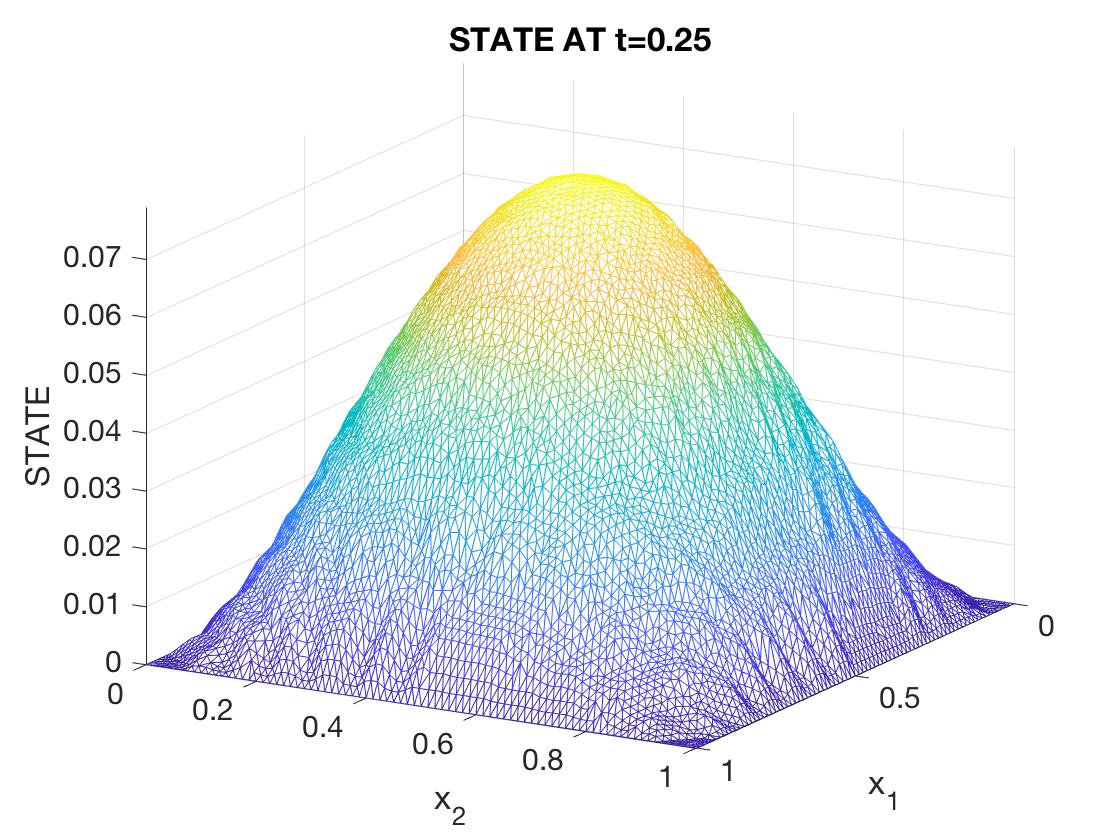}\quad
\includegraphics[height=48mm, width = 62mm]{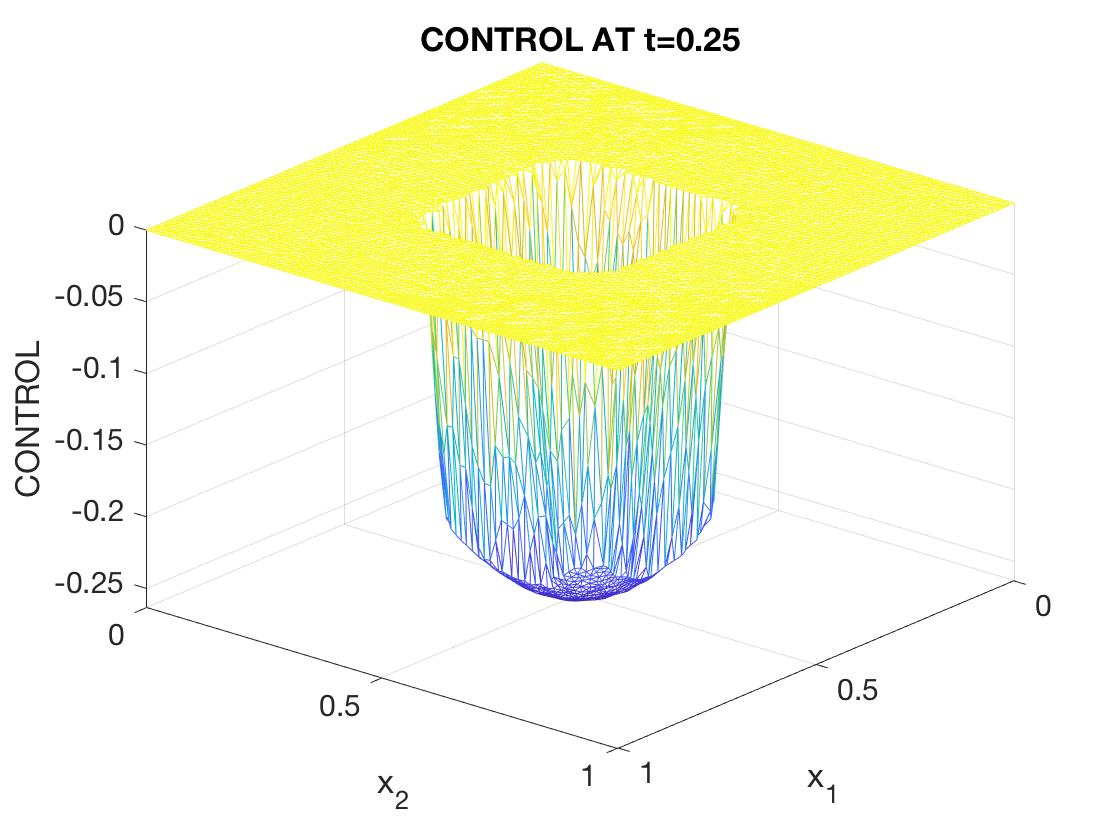}
\caption{Test 2.1 -- Cuts of the computed state and control at $t = 0.25$.}
\label{Fig_2.1-4}
\end{figure}

\begin{figure}[htbp]
\centering
\includegraphics[height=48mm, width = 62mm]{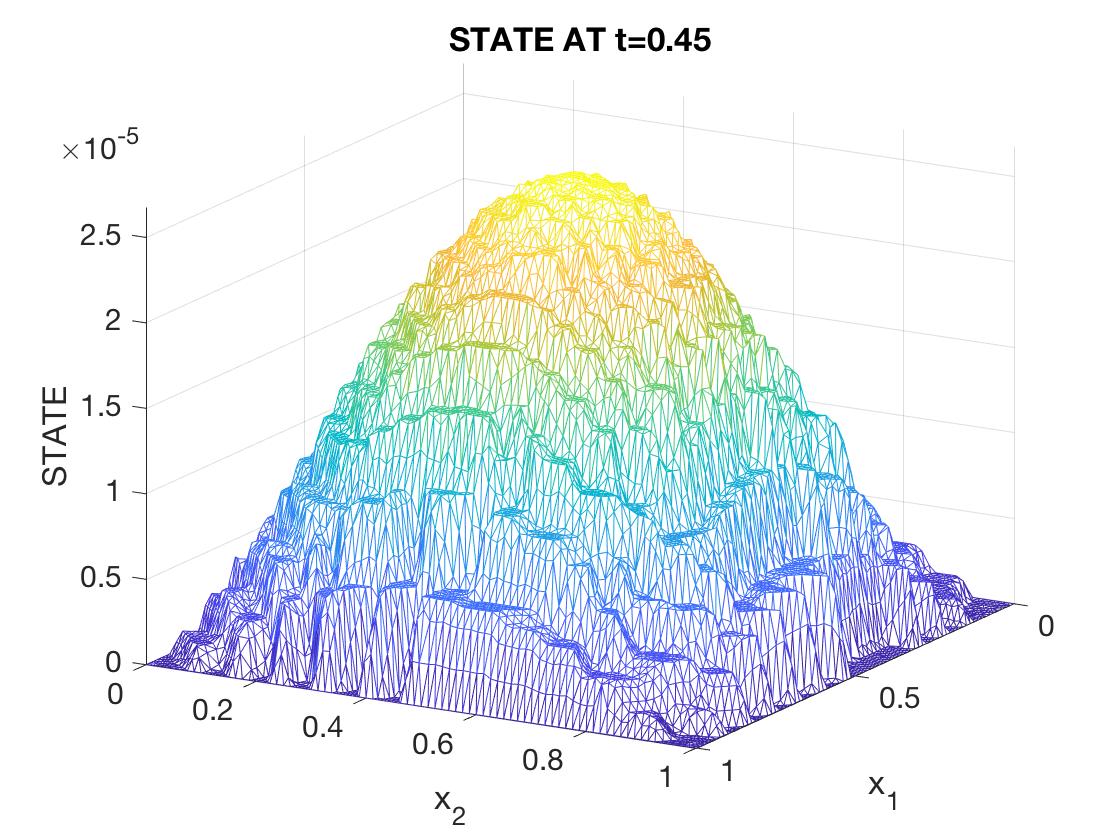}\quad
\includegraphics[height=48mm, width = 62mm]{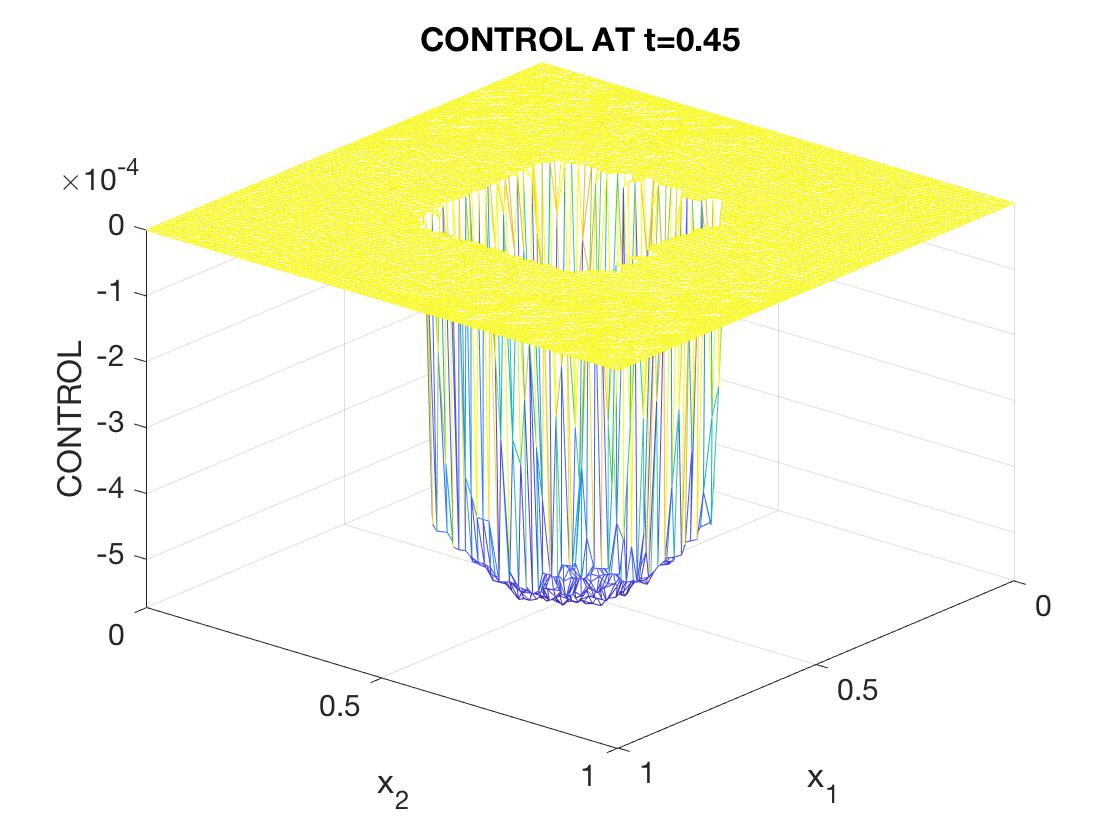}
\caption{Test 2.1 -- Cuts of the computed state and control at $t = 0.45$.}
\label{Fig_2.1-5}
\end{figure}

\begin{figure}[htbp]
\centering
\includegraphics[height=48mm, width = 62mm]{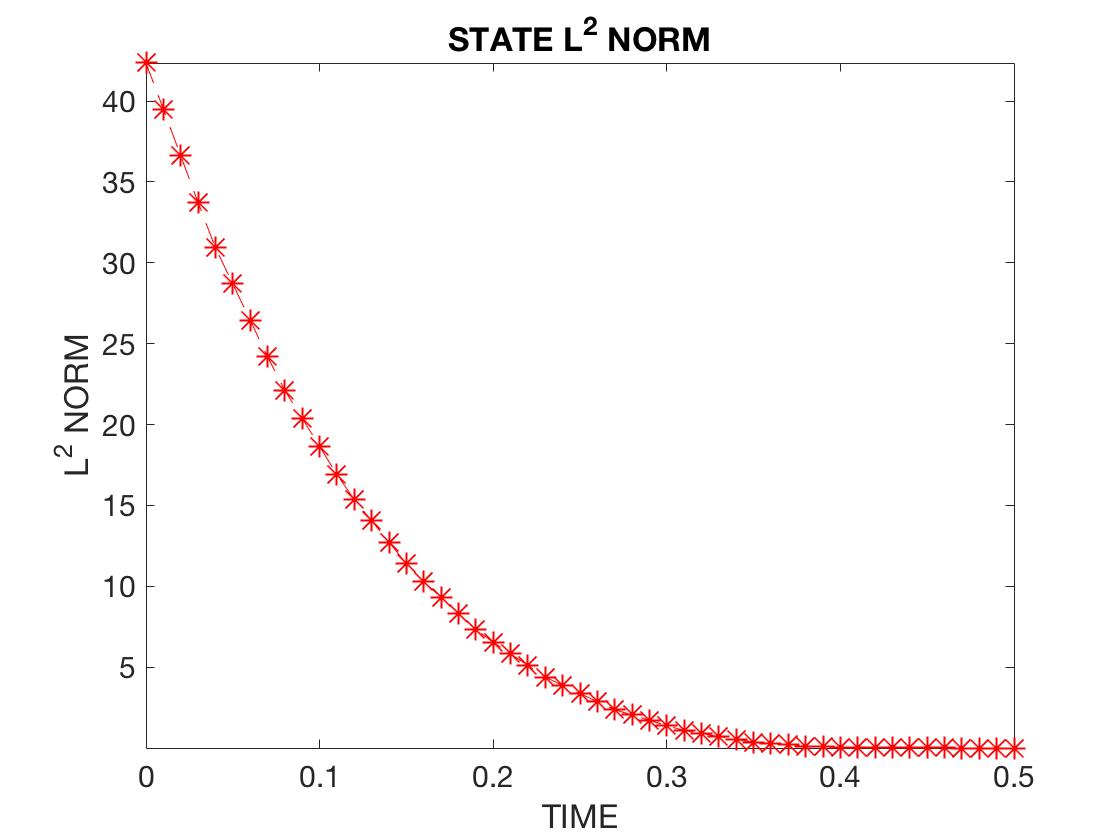}\quad
\includegraphics[height=48mm, width = 62mm]{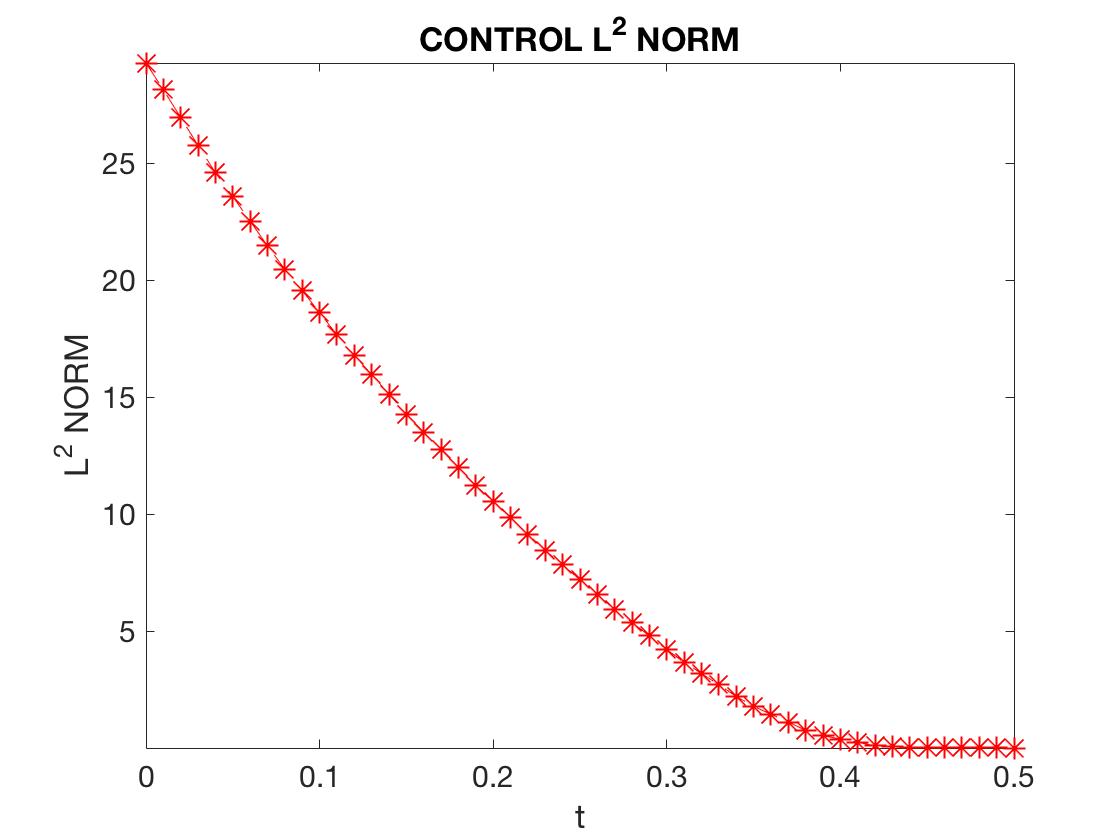}
\caption{Test 2.1 -- Evolution in time of the $L^2$ norms of the state and the control.}
\label{Fig_2.1-6}
\end{figure}

   Again, the solution has been computed for several mesh sizes $h=(\Delta x_1,\Delta x_2,\Delta t)$.
   We show the cylinder~$Q$ and a mesh in~Figure~\ref{Fig_2.1-1}.
   Several views of the computed state and control are displayed in~Figures~\ref{Fig_2.1-3}--\ref{Fig_2.1-5}.
   Also, the evolution in time of the $L^2$ norms of the state and the control is shown in~Figure~\ref{Fig_2.1-6}.

   Starting again from $(y^0, v^0)= H'(0,0)^{-1} (0,y_0)$, we have used \textbf{ALG~1} to compute the solution to the null control problem.
   The evolution of the relative errors of the computed solutions for several mesh sizes is given in Table~\ref{Table_2.1-1};
   again, we have assumed that the ``exact'' solution corresponds to the finest mesh.

   The number of iterates needed to get convergence is indicated in Table~\ref{Table_2.1-2} for several mesh sizes.
   Finally, we give in Table~\ref{Table_2.1-3} this number for $y_0 \equiv y_{00} \sin(\pi x_1) \sin(2\pi x_2)$ and several choices of~$y_{00}$ when the finest mesh is used.

\begin{table}[htbp]
\centering
\caption{Test 2.1 -- The control and state relative errors for various $h = (\Delta x_1,\Delta x_2,\Delta t)$.}
\label{Table_2.1-1}
\begin{tabular}{ccccc}
\toprule
$h = (\Delta x_1,\Delta x_2,\Delta t)$  & Points & Elements & Control rel error & State rel error \\
\midrule
$(0.250,0.250,0.056)$  &  $580$ & $2646$ & $4.5015 \times 10^{-1}$ & $2.7034 \times 10^{-1}$ \\
$(0.111,0.111,0.056)$  &  $2590$ & $12960$ & $9.0513 \times 10^{-2}$ & $6.9357 \times 10^{-2}$ \\
$(0.077,0.077,0.038)$  &  $5096$ & $26286$ & $4.5896 \times 10^{-2}$ & $2.7980 \times 10^{-2}$ \\
$(0.053,0.053,0.026)$  &  $14900$ & $80484$ & $9.6015 \times 10^{-2}$ & $1.0923 \times 10^{-2}$ \\
$(0.048,0.048,0.024)$  &  $20482$ & $111888$ & $5.4722 \times 10^{-3}$ & $6.1146 \times 10^{-3}$ \\
$(0.040,0.040,0.020)$  &  $33384$ & $185950$ & $1.2821 \times 10^{-5}$ & $6.5718 \times 10^{-5}$ \\
$(0.034,0.034,0.017)$  &  $50290$ & $281010$ &  --- & --- \\
\bottomrule
\end{tabular}
\end{table}

\begin{table}[htbp]
\centering
\caption{Tests 2.1 -- Iterates needed to get~\eqref{stopping} (convergence) for various $\kappa$ and~$h = (\Delta x_1,\Delta x_2,\Delta t)$.}
\label{Table_2.1-2}
\begin{tabular}{ccccc}
\toprule
$h = (\Delta x_1,\Delta x_2,\Delta t)$  & Points & Elements & Iter.\ for $\kappa=10^{-5}$ & Iter.\ for $\kappa=10^{-9}$ \\
\midrule
$(0.250,0.250,0.056)$  &  $580$ & $2646$ & 15 & 24 \\
$(0.111,0.111,0.056)$  &  $2590$ & $12960$ & 17 & 24 \\
$(0.077,0.077,0.038)$  &  $5096$ & $26286$ & 20 & 26 \\
$(0.053,0.053,0.026)$  &  $14900$ & $80484$ & 23 & 28 \\
$(0.048,0.048,0.024)$  &  $20482$ & $111888$ & 30 & 45 \\
$(0.040,0.040,0.020)$  &  $33384$ & $185950$ & 41 & 69 \\
$(0.034,0.034,0.017)$  &  $50290$ & $281010$ & 67 & 123 \\
\bottomrule
\end{tabular}
\end{table}

\begin{table}[htbp]
\centering
\caption{Tests 2.1 -- Iterates needed to get~\eqref{stopping} (convergence) for various $\kappa$ and initial data.}
\label{Table_2.1-3}
\begin{tabular}{ccc}
\toprule
$y_{00}$ & Iter.\ for $\kappa=10^{-5}$ & Iter.\ for $\kappa=10^{-9}$ \\
\midrule
$ 0.25 $  &  $11$ & $22$ \\
$ 0.50 $  &  $25$ & $50$ \\
$ 0.75 $  &  $41$ & $69$ \\
$ 1.00 $  &  $123$ & $170$ \\
$ 1.25 $  &  $234$ & $350$ \\
\bottomrule
\end{tabular}
\end{table}

\subsection{Test 2.2}

   In this last Test, we take the same data as in Test~2.1, except for the diffusion coefficient.
   Now,

\begin{itemize}
\item $a(s) \equiv 1 + 0.7 |s|^2$.
\end{itemize}

   We can obtain similar numerical results by applying \textbf{ALG~1} in combination with a mixed finite element approximation.
   The computations lead to Tables~\ref{Table_2.2-1} and~\ref{Table_2.2-2}.

\begin{table}[htbp]
\centering
\caption{Test 2.2 -- The control and state relative errors.}
\label{Table_2.2-1}
\begin{tabular}{ccccc}
\toprule
$h = (\Delta x_1,\Delta x_2,\Delta t)$  & Points & Elements & Control rel error & State rel error \\
\midrule
$(0.250,0.250,0.056)$  &  $580$ & $2646$ & $3.0068 \times 10^{-1}$ & $5.6655 \times 10^{-1}$ \\
$(0.111,0.111,0.056)$  &  $2590$ & $12960$ & $1.7731 \times 10^{-1}$ & $9.8043 \times 10^{-2}$ \\
$(0.077,0.077,0.038)$  &  $5096$ & $26286$ & $5.4722 \times 10^{-2}$ & $3.7850 \times 10^{-2}$ \\
$(0.053,0.053,0.026)$  &  $14900$ & $80484$ & $1.5887 \times 10^{-2}$ & $9.0923 \times 10^{-3}$ \\
$(0.048,0.048,0.024)$  &  $20482$ & $111888$ & $5.2652 \times 10^{-3}$ & $2.1465 \times 10^{-3}$ \\
$(0.040,0.040,0.020)$  &  $33384$ & $185950$ & $1.2821 \times 10^{-4}$ & $1.0804 \times 10^{-4}$ \\
$(0.034,0.034,0.017)$  &  $50290$ & $281010$ &  --- & --- \\
\bottomrule
\end{tabular}
\end{table}

\begin{table}[htbp]
\centering
\caption{Tests 2.2 -- Iterates needed to get~\eqref{stopping} (convergence) for various $\kappa$ and~$h = (\Delta x_1,\Delta x_2,\Delta t)$.}
\label{Table_2.2-2}
\begin{tabular}{ccccc}
\toprule
$h = (\Delta x_1,\Delta x_2,\Delta t)$  & Points & Elements & Iter.\ for $\kappa=10^{-5}$ & Iter.\ for $\kappa=10^{-9}$ \\
\midrule
$(0.250,0.250,0.056)$  &  $580$ & $2646$ & 25 & 44 \\
$(0.111,0.111,0.056)$  &  $2590$ & $12960$ & 28 & 44 \\
$(0.077,0.077,0.038)$  &  $5096$ & $26286$ & 31 & 46 \\
$(0.053,0.053,0.026)$  &  $14900$ & $80484$ & 53 & 78 \\
$(0.048,0.048,0.024)$  &  $20482$ & $111888$ & 60 & 85 \\
$(0.040,0.040,0.020)$  &  $33384$ & $185950$ & 65 & 89 \\
$(0.034,0.034,0.017)$  &  $50290$ & $281010$ & 97 & 150 \\
\bottomrule
\end{tabular}
\end{table}

\FloatBarrier

\section{Conclusions, further remarks and open questions}\label{Sec_Comments}

   We have established the local null controllability property of the quasi-linear parabolic system~\eqref{sistema1fase}, at any given positive time we prescribe, both with distributed and boundary controls.
   We found that a general method, which is quite standard nowadays ---
   see~\cite{CCLM, fernandez2015theoretical, fernandez2021theoretical, limaco2016null}, to name a few --- is also successful in the present context.

    In what concerns theoretical control, the main novelty in this paper is the proof of new nontrivial regularity results of the control (and the state), as well as the boot-strapping and verification arguments --- the aspects that make it possible to use these techniques.

    Indeed, the local inversion argument needed to prove Theorem~\ref{Maintheorem} only works after introducing the spaces~$Y$ and~$Z$ in~\eqref{3.1a}--\eqref{3.1aa} and checking that the mapping $H$ in~\eqref{defH} is well-defined and~$C^1$ and~$H'(0,0)$ is onto.

   By inspection, we see that the proof in Section~\ref{SecLocNC1Phase} can be adapted to deal with the system
   $$
\begin{cases}
y_t - \nabla \cdot [A(y,\nabla y) D^2 y] = f(y,\nabla y) + \chi_\omega v\ &\text{ in } Q, \\
y = 0\ &\text{ on } \Sigma, \\
y|_{t=0} = y_0\ &\text{ in } \Om,
\end{cases}
   $$
where $A : \mathbb{R} \times \mathbb{R}^d \mapsto \mathbb{R}^{d \times d}$ and~$f : \mathbb{R} \times \mathbb{R}^d \mapsto \mathbb{R}$ satisfy the following assumptions:
\begin{itemize}
   \item[$\textbf{H1}^{\prime}$] $A$ is of class $C^4$ and~$f$ is of class $C^3$.
   \item[$\textbf{H2}^{\prime}$] There exist $C$, $a_0 >0$ and~$r,s,r^\prime,s^\prime \geq 1$ such that
   $$
   \begin{cases}
   \xi^{T} A(u,p) \xi \geq a_0|\xi|^2, \\
   |\partial_1^\alpha D_2^{\beta} a(u,p)| \leq C\left(1 + |u|^{(s-|\alpha|)^+} + |p|^{(r-|\beta|)^+} \right) \ \text{ and} \\
   |\partial_1^{\gamma}D_2^{\sigma} f(u,p) | \leq C\left(1 + |u|^{(s^\prime-|\gamma|)^+} + |p|^{(r^\prime-|\sigma|)^+} \right),
   \end{cases}
   $$
for all $(u,p,\xi) \in \mathbb{R}\times \mathbb{R}^d\times \mathbb{R}^d$, all multi-indices $\beta,\sigma$ subject to the constraints $|\beta|\leq 4$ and~$|\sigma|\leq 3$, and all integers $\alpha, \gamma$ satisfying $|\alpha|\leq 4$ and~$|\gamma|\leq 3$.
   Above, $\partial_1$ denotes the partial derivative with respect to the one  dimensional variable $u$, while $D_2$ refers to derivatives with respect to the coordinates of $p \in \mathbb{R}^d$.
\end{itemize}

   Furthermore, we can also consider systems in non-divergence form, namely,
   $$
\begin{cases}
y_t - A(y,\nabla y) : D^2 y = f(y,\nabla y) + \chi_\omega v\ \text{ in } Q, \\
y = 0\ \text{ on } \Sigma, \\
y|_{t=0} = y_0\ \text{ in } \Om,
\end{cases}
   $$
where the assumptions on $f$ remain the same but less regularity on $A$ is required.
   More precisely, it suffices to assume that $A$ is of class $C^3$ and satisfies power law limitations, like in~\textbf{H2} or~$\textbf{H2}^\prime$, but now up to derivatives of order three.

   Another natural question is whether Theorem~\ref{Maintheorem} still holds for similar systems with PDE's of the form
   \begin{equation}\label{non-hom}
y_t -\nabla \cdot (a(t,x;\nabla y)\nabla{y})= v \tilde{1}_\omega,
   \end{equation}
that is, with a nonlinear diffusion coefficient non-homogeneous in space or time.
   The answer is affirmative as long as we suppose enough regularity for $a : \overline{Q} \times \mathbb{R}^d \mapsto \mathbb{R}$.
   It is also worthwhile to note that the local exact controllability to the (smooth) trajectories is an interesting additional topic.
   After a change of variable, this can be reduced to the null controllability of a system where the PDE is similar to~\eqref{non-hom}
   (see some related comments in~Section~\ref{Sec_Comments}).
   For brevity, we will not delve into the technical aspects of this extension.

   It is possible to demonstrate local null controllability properties for multi-phase quasi-linear parabolic systems, assuming uniform ellipticity of the diffusion coefficients and that all phases are to be controlled.
   For two-phase systems, it is also straightforward to treat the case of zero-order coupling with one distributed internal control force.
   For more details and some related results, see for instance~\cite{Barbu-FB, Krstic}.

   The techniques developed here may serve to address more general models of chemotaxis, as in~\cite{hillen2009user}, possibly leading to results similar to those in~\cite{chaves2015uniform,chaves2017controllability}.
   In~\cite{de2019local}, the authors showed that the present method also works in addressing the local null controllability of a generalization of the Navier-Stokes equations.

   Regarding the limitations of the present method, the drawbacks are the dimensional constraint $d \leq 3$ and the higher regularity demanded to the initial data.

   Figures~\ref{Fig_1.1-1}, \ref{Fig_1.2-1} and~\ref{Fig_1.3-3} are similar to Figures~3 and~4 in~\cite{FCMS} (that correspond to the solution of a null controllability problem for the linear 1D heat equation), Figure~7 in~\cite{FCMunch} (corresponding to a semilinear heat equation) and Figures~4, 5, 7 and~8 in~\cite{fernandez2000null} (corresponding to a quasi-linear parabolic equation where the diffusion coefficient depends on the state and not its gradient).
   Similar comments can be made for the results of Test~2.1.
   Thus, Figures~\ref{Fig_2.1-3} to~\ref{Fig_2.1-5} are similar to Figures~2 and~3 in~\cite{CHCV} (dealing with linear 2D heat equations) and (for instance) Figure~7 in~\cite{FCMunchNewton} (corresponding to a quasi-linear parabolic equation).

\begin{figure}[htbp]
\centering
\includegraphics[height=48mm, width = 62mm]{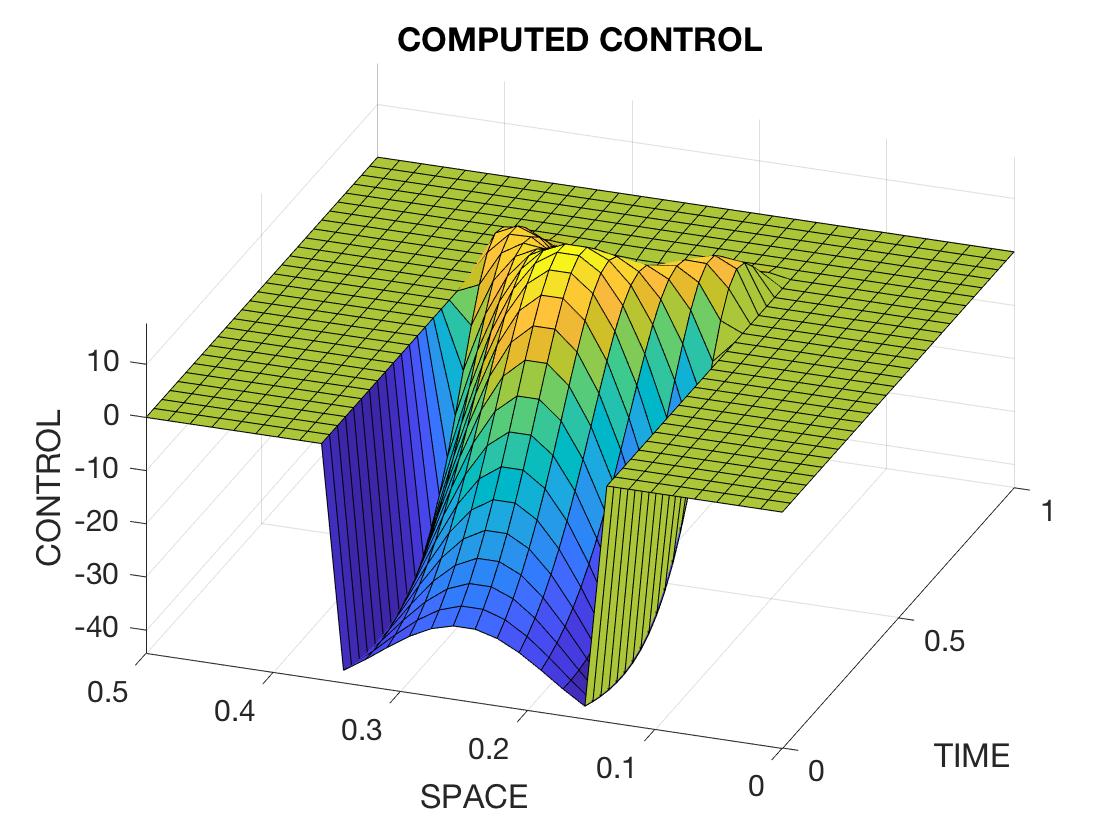}\quad
\includegraphics[height=48mm, width = 62mm]{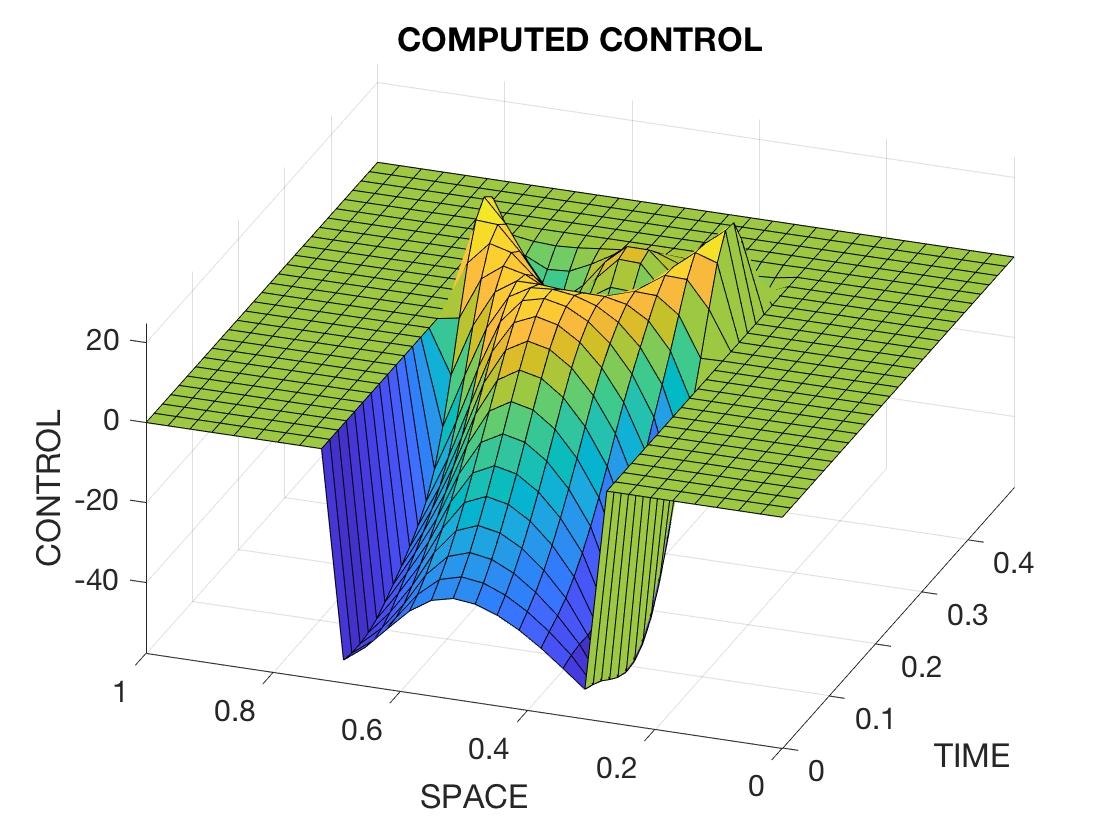}
\caption{Additional Test -- The computed controls corresponding the linear case ({\it left}) and the quasi-linear case ({\it right}).}
\label{Fig_Add-1}
\end{figure}

\begin{figure}[htbp]
\centering
\includegraphics[height=48mm, width = 62mm]{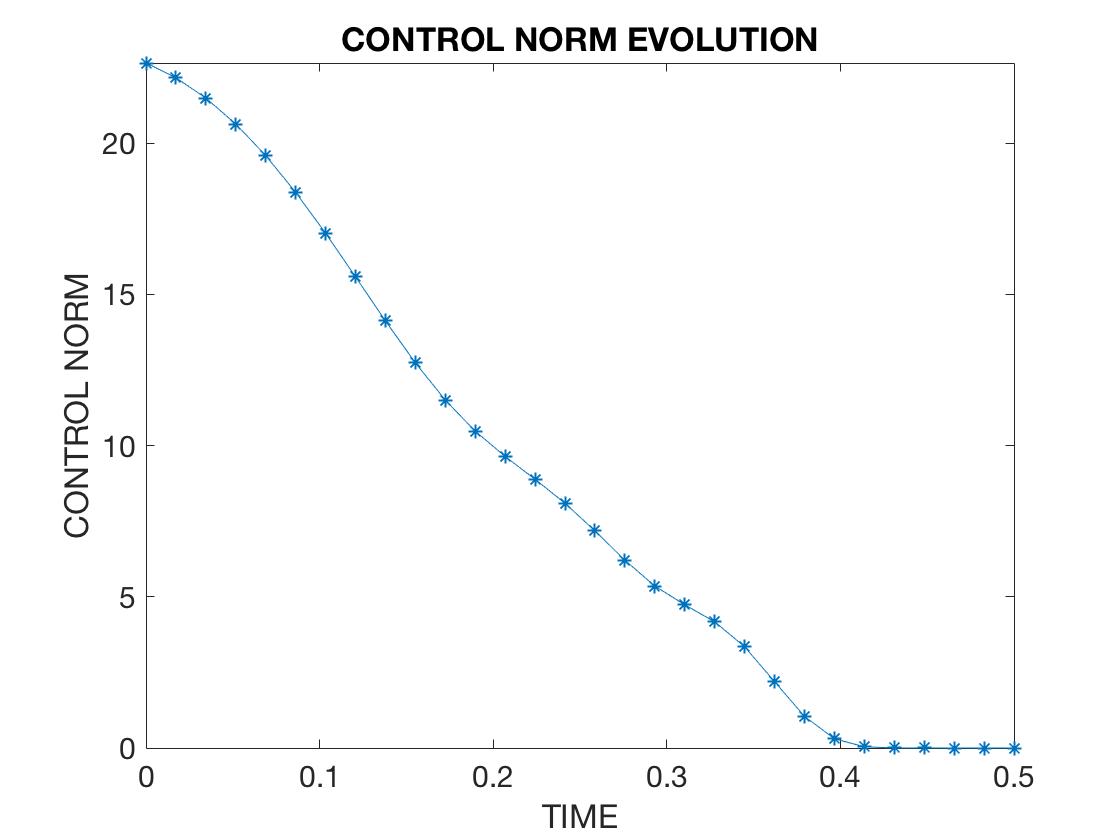}\quad
\includegraphics[height=48mm, width = 62mm]{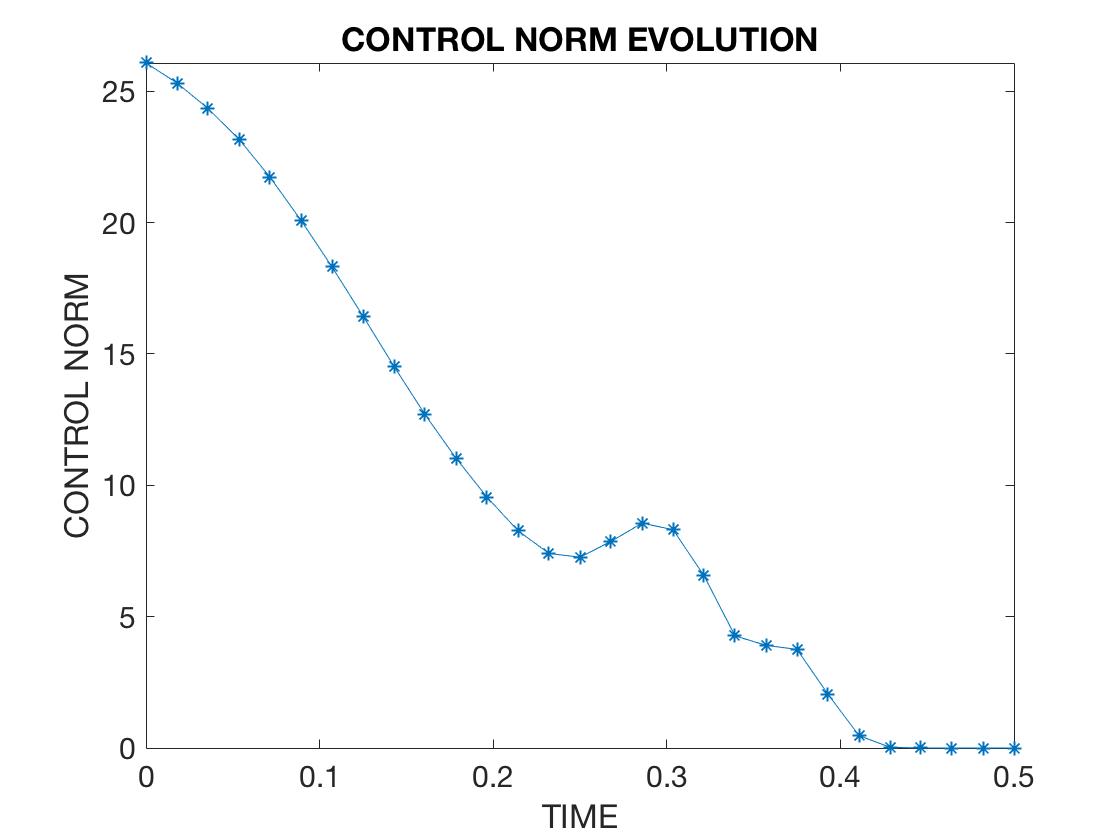}
\caption{Additional Test -- Evolution in time of the $L^2$ norms of the computed controls corresponding the linear case ({\it left}) and the quasi-linear case ({\it right}).}
\label{Fig_Add-2}
\end{figure}

\begin{figure}[htbp]
\centering
\includegraphics[height=48mm, width = 62mm]{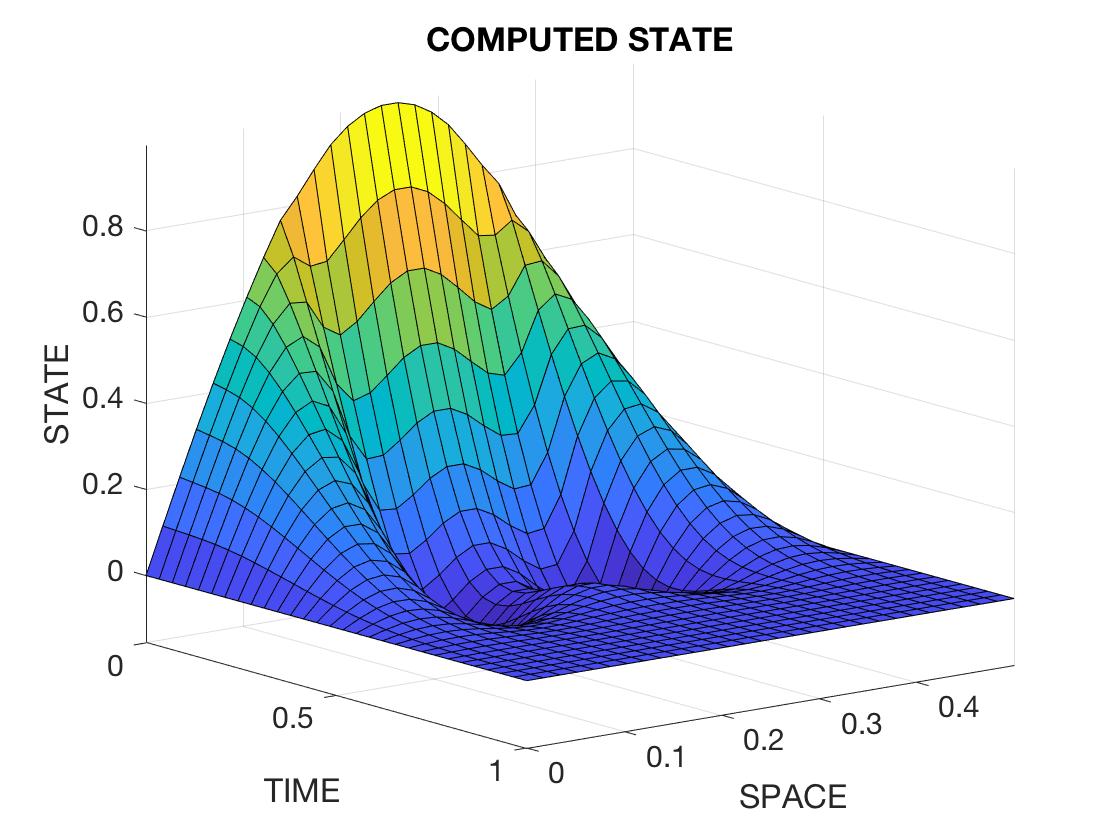}\quad
\includegraphics[height=48mm, width = 62mm]{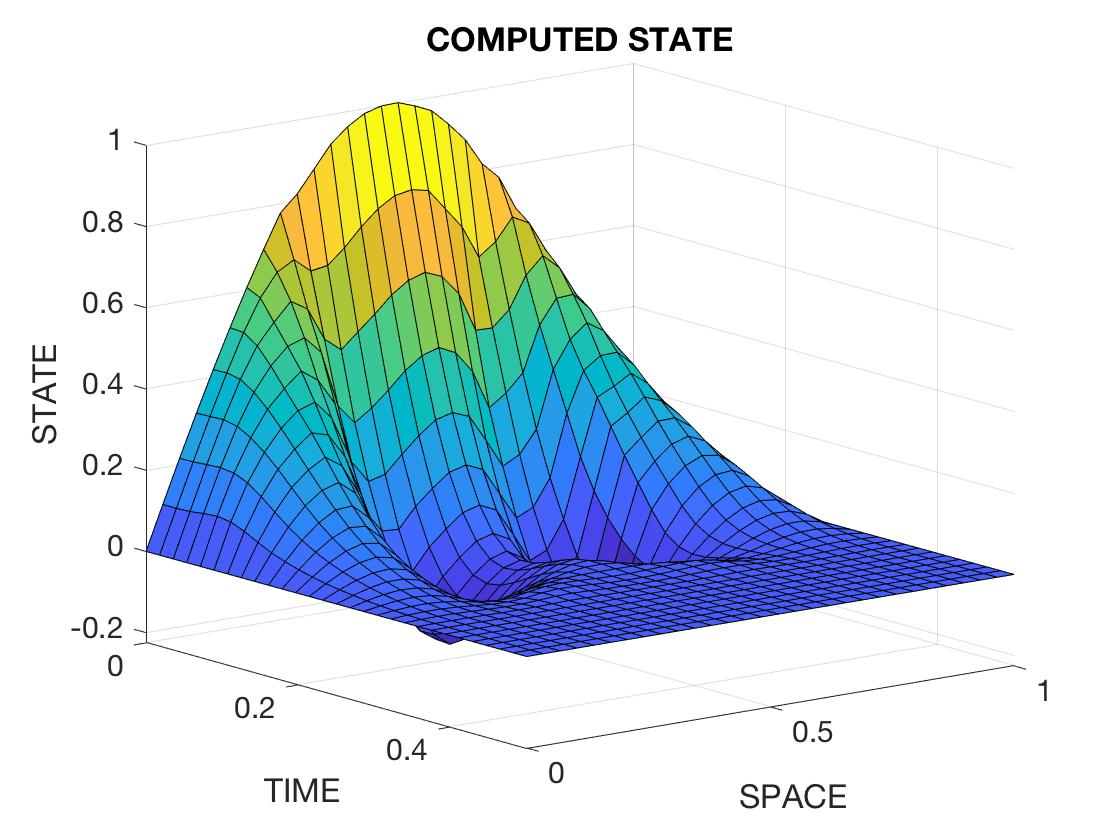}
\caption{Additional Test -- The computed states corresponding the linear case ({\it left}) and the quasi-linear case ({\it right}).}
\label{Fig_Add-3}
\end{figure}

\begin{figure}[htbp]
\centering
\includegraphics[height=48mm, width = 62mm]{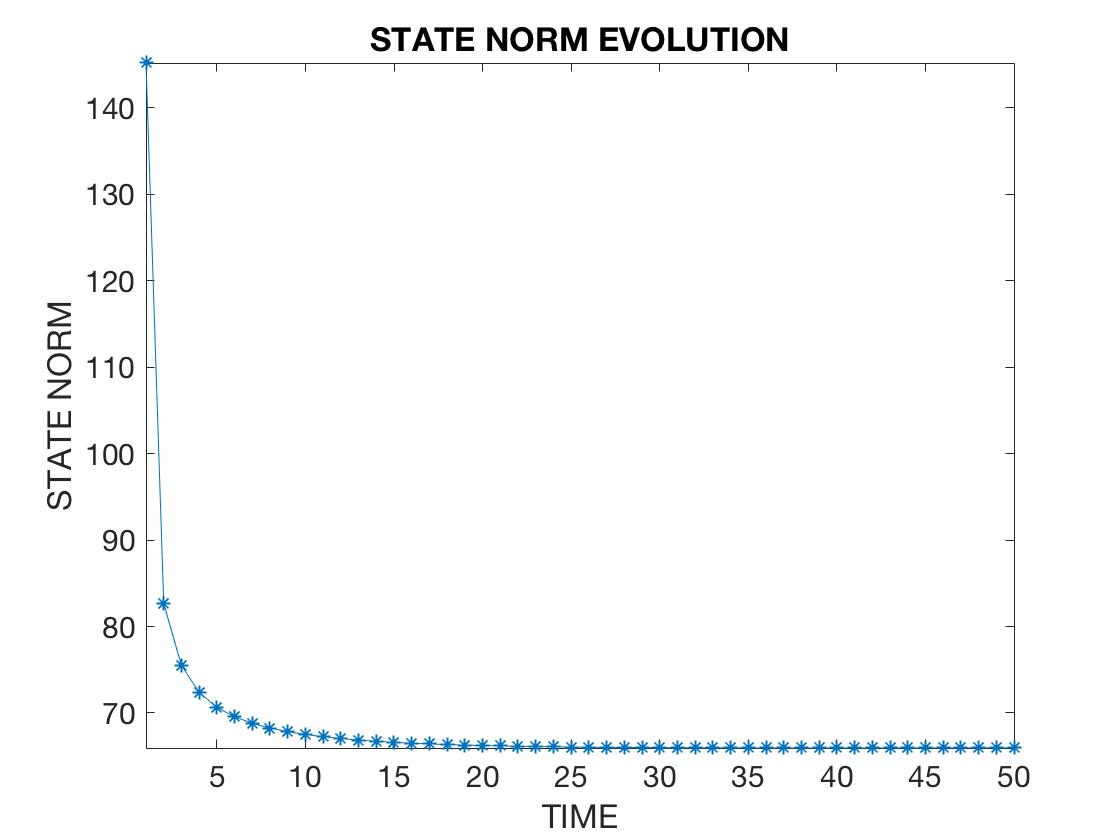}\quad
\includegraphics[height=48mm, width = 62mm]{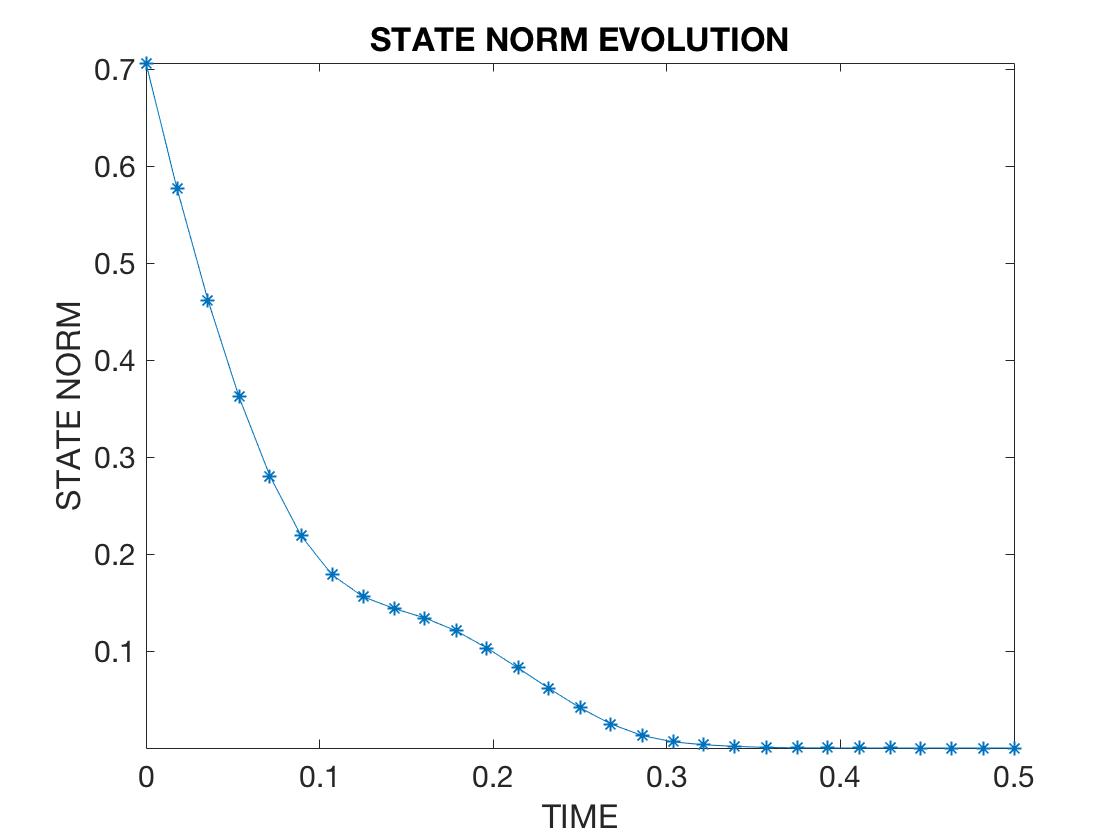}
\caption{Additional Test -- Evolution in time of the $L^2$ norms of the computed states corresponding the linear case ({\it left}) and the quasi-linear case ({\it right}).}
\label{Fig_Add-4}
\end{figure}

   In both cases, the effect of nonlinearity is observed.
   To make this assertion more convincing, we have depicted in Figures~\ref{Fig_Add-1}--\ref{Fig_Add-4} the computed controls corresponding to the solution of~\eqref{sistema1fase} with

\begin{itemize}
\item $d=1$, $\Om = ]0,1[$, $\omega = ]0.3, 0.7[$, $T = 0.5$.
\item $y_0(x) \equiv \sin(\pi x)$.
\item $a(s) \equiv 0.1$ (linear case) and $a(s) \equiv 0.1 + 0.05 \sin(s)$ (quasi-linear case).
\end{itemize}

   Note in particular that the state norms displayed in Figure~\ref{Fig_Add-4} evolve in a clearly different way in the linear and quasi-linear cases.

   Finally, recall that there are many open problems concerning the null control of semilinear and nonlinear systems.
   In particular, as already mentioned at the beginning of the paper, it is unknown whether the global controllability property holds for~\eqref{sistema1fase}.

\FloatBarrier

\appendix

\section{Proof of Corollary~\ref{RegularityForLinear1PhaseProblem}}\label{Sec-App-A}

   For the proof of the assertions in~\textit{(a)} and~\textit{(b)}, see~\cite[Theorem~2.1]{fernandez2021theoretical}.
   The claims in~\textit{(c)} and~\textit{(d)}, which we are going to prove, are new. 
   
   With the notations introduced after the statement of~Theorem~\ref{ControlOfLinearProblem}, we set
   \begin{align*}
   \rho_{k}\hat{v}_t &= 2s\overline{\alpha} e^{-s\overline{\alpha}/\ell} \ell^{-5 + k/2}l^\prime p - 3 e^{-s\overline{\alpha}/\ell} \ell^{-4 + k/} \ell^\prime \varphi + e^{-s\overline{\alpha}/\ell} \ell^{-3 + k/2} p_t \\
   &=: A_1 + A_2 + A_3,
   \end{align*} 
where $k \geq 1$ is an integer that will be fixed below.
   Therefore, we have $L^*(\rho_k \hat{v}_t) = L^*A_1 + L^*A_2 + L^*A_3$.
   
   It is clear that 
   $$
|L^*A_1 + L^*A_2| \leq C \left( \rho_{10-k}^{-1}|L^*p| + \rho_{14-k}^{-1}|p| + \rho_{12-k}^{-1}|\nabla p| \right),
   $$
   $$
|L^*A_3| \leq C\left( \rho_{6-k}^{-1}|L^*p_t| + \rho_{10-k}^{-1}|p_t| + \rho_{8-k}^{-1}|\nabla p_t| \right).
   $$
   For $k \geq 11$, it follows that $L^*A_1 + L^*A_2 \in L^2(Q)$, with
   $$
\|L^*A_1 + L^*A_2\|_{L^2(Q)}^2 \leq C\pi(p,p) \leq C\left( \|y_0\|^2 + \jjntQ \rho_3^2 |f|^2\  \right).
   $$
   On the other hand, $L^*p_t = (e^{2s\overline{\alpha}/\ell} y)_t = e^{2s\overline{\alpha}/\ell} y_t - 2s\overline{\alpha} 
e^{2s\overline{\alpha}/\ell} \ell^{-2} \ell^\prime y$ and, consequently,
   $$
|L^*p_t| \leq C \left( \rho^2 |y_t| + \rho_{-2}^{2}|y| \right).
   $$
   From the estimate~\eqref{regularity1}, we can conclude that $\rho_{-7}^{-1}L^*p_t \in L^2(Q)$.
   Moreover, from the Carleman inequality~\eqref{CI-new}, we have~$\rho_{-1}^{-1}p_t \in L^2(Q)$.
   Finally, we observe that $\rho_{-7}^{-1}L^*p_t \in L^2(Q)$ and~$\rho_{-1}^{-1}p_t \in L^2(Q)$ imply $\rho_{-7}^{-1}\nabla p_t \in \left[L^2(Q)\right]^d$.
   Hence, for $k \geq 15$, it is true that $L^*A_3 \in L^2(Q)$ and
   $$
\|L^*A_3\|_{L^2(Q)}^2 \leq C\left[ \pi(p,p) + \|\rho_7 y_t\|_{L^2(Q)}^2 + \|\rho y\|_{L^2(Q)}^2 \right].
   $$
   Recalling that $\rho_{15} \widehat{v} = 0$ on~$\Sigma$ and~$(\rho_{15}\hat{v})(T\,,\cdot) = 0$ in~$\Om$, we conclude that $v = \chi\hat{v}|_{]0,T[ \times \omega}$ satisfies $(\rho_{15}v_t)_t \in L^2(]0,T[ \times \omega)$ and~$\rho_{15} v_t \in L^2(0,T;H^2(\omega))$ and the corresponding estimates are fulfilled.

   In a similar way, we can show that $L^*(\rho_{15}\Delta \hat{v}) \in L^2(Q)$.
   This, together with the identities $\rho_{15} \Delta\widehat{v} = 0$ on~$\Sigma$ and~$(\rho_{15}\Delta\hat{v})(T\,,\cdot) = 0$ in~$\Om$, give $(\rho_{15}\Delta v)_t \in L^2(]0,T[ \times \omega)$ (which we already knew), $\rho_{15}\Delta v \in L^2(0,T;H^4(\omega))$ and~(\ref{regularity3}).

   Let us establish the estimates \textit{(d).} 
   
   Henceforth, we fix $n\geq 1$, and we let $v_n,\, f_n$ and $y_{0n}$ denote the projections of $v,\, f$ and $y_0$ on the space spanned by the first $n$ eigenfunctions of the Dirichlet Laplacian in~$\Om$.
   Then, we let $y_n$ denote the corresponding solution to~\eqref{Sistema1FaseLinear}.
   
   For the sake of brevity, we drop in the next computations the indices $n$ in $v_n,f_n, y_{0n}$ and $y_n$.
   The constants $C$ arising in the estimates below are independent of $n$.
   
   We divide the proof into four parts.

\textsc{Part~1:} We differentiate system (\ref{Sistema1FaseLinear}) twice with respect to time.
   It follows that
   \begin{equation} \label{LinearTwoTimeDerivs}
   y_{ttt} - \sigma\Delta y_{tt} = \chi_\omega v_{tt} + f_{tt}.
   \end{equation}
   By multiplying the above equation by $\rho_{15}^2 y_{tt}$ and integrating in space over $\Om$, after some integrations by parts, we obtain:
   $$
   \begin{array}{l} \displaystyle
\frac{1}{2}\frac{d}{dt}\left(\jnt_\Om \rho_{15}^2 |y_{tt}|^2 \,dx \right) + \sigma \jnt_\Om \rho_{15}^2 |\nabla y_{tt}|^2\,dx
= \jnt_\Om \rho_{15}^2 v_{tt}y_{tt}\,dx + \jnt_\Om \rho_{15}^2 f_{tt}y_{tt}\,dx
\\ \noalign{\smallskip} \displaystyle
\hspace{1.5cm} \ + \frac{1}{2}\jnt_\Om (\rho_{15}^2)_t |y_{tt}|^2\,dx - \sigma \jnt_\Om \nabla(\rho_{15}^2)\cdot\nabla y_{tt} y_{tt} \,dx 
\\ \noalign{\smallskip} \displaystyle
\hspace{1.3cm} \ \leq C\left(\jnt_\Om \rho_{15}^2 |v_{tt}|^2\,dx + \jnt_\Om  \rho_{17}^2 |f_{tt}|^2\,dx
+ \jnt_\Om \rho_{13}^2 |y_{tt}|^2\,dx \right) + \frac{\sigma}{2}\jnt_\Om \rho_{15}^2 |\nabla y_{tt}|^2\,dx.
\end{array}   
   $$
   
   Note that the identity $y_{tt}(0\,,\cdot) = \Delta^2 y_0 + \chi_\omega \Delta v(0\,,\cdot) + \Delta h(0\,,\cdot) + \chi_\omega v_t(0\,,\cdot) + f_t(0\,,\cdot)$ yields $\|y_{tt}|_{t=0}\| \leq S(y_0,v,f)$, where $S(y_0,v,f)$ denotes the right hand side of~\eqref{regularity3}.
   Thus, using these relations and the estimate~\eqref{regularity3}, we deduce for the Galerkin approximations and then, by standard limiting arguments, for the actual solution that
   \begin{equation} \label{2.15p}
\sup_{\left[0,T\right]}\left( \jnt_\Om \rho_{15}^2 |y_{tt}|^2 \,dx \right) + \jjntQ \rho_{15}^2 |\nabla y_{tt}|^2 
\leq S(y_0,v,f).
   \end{equation}

\textsc{Part~2:} Now, multiplying (\ref{LinearTwoTimeDerivs}) by~$\rho_{17}^2(y_{ttt}-\sigma\Delta y_{tt})$, integrating in space and performing integrations by parts, we find the following:
   \begin{equation} \label{DontKnowHowToCallThis}
\begin{array}{l} \displaystyle
   \jnt_\Om \rho_{17}^2\left[ |y_{ttt}|^2 + \sigma^2  (\Delta y_{tt})^2 \right]dx + \sigma \frac{d}{dt}\left(\jnt_\Om \rho_{17}^2 |\nabla y_{tt}|^2\,dx \right) 
   \\ \noalign{\smallskip} \displaystyle
   \hspace{0.7cm} \ = \jnt_\Om \rho_{17}^2 v_{tt}(y_{ttt} - \sigma\Delta y_{tt})\,dx + \jnt_\Om  \rho_{17}^2 f_{tt}(y_{ttt}-\sigma\Delta y_{tt}) \,dx 
   \\ \noalign{\smallskip} \displaystyle
   \hspace{0.7cm} \ + \sigma \jnt_\Om  (\rho_{17}^2)_t|\nabla y_{tt}|^2\,dx - 2\sigma \jnt_\Om y_{ttt} \nabla(\rho_{17}^2)\nabla y_{tt} \,dx
   \\ \noalign{\smallskip} \displaystyle
   \hspace{0.7cm} \ \leq  C\left(\jnt_\Om \rho_{17}^2 |v_{tt}|^2\,dx + \jnt_\Om \rho_{17}^2 |h_{tt}|^2\,dx + \jnt_\Om \rho_{15}^2 |y_{tt}|^2\,dx \right) + \frac{1}{2} \jnt_\Om  \rho_{17}^2\left[ y_{ttt}^2 + \sigma^2  (\Delta y_{tt})^2\right]\,dx.
\end{array}   
   \end{equation}
   From~\eqref{2.15p}, the identity $\nabla y_{tt}(0\,,\cdot) = \nabla \Delta^2 y_0 + \chi_\omega \nabla\Delta v(0\,,\cdot) + \nabla\Delta h(0\,,\cdot) + \nabla v_t(0\,,\cdot) + \nabla f_t(0\,,\cdot)$ and the fact that $\rho_{15}v \in L^\infty(0,T;H^3(\omega)),\ \rho_{15}v_t \in L^\infty(0,T;H^1_0(\omega))$ and
   $$
\sup_{\left[0,T\right]}\left\{ \|\rho_{15}v\|_{H^3} + \|\rho_{15}v_t\|_{H^1}\right\} \leq S(y_0,v,f),
   $$
after integration in time of~(\ref{DontKnowHowToCallThis}), we deduce for the Galerkin approximates, that 
   \begin{equation} \label{2.16p}
\jjntQ \rho_{17}^2\left[ |y_{ttt}|^2 + (\Delta y_{tt})^2 \right]  
+ \sup_{\left[0,T\right]}\left[ \jnt_\Om \rho_{17}^2 |\nabla y_{tt}|^2\,dx\right](t)
\leq C S(y_0,v,f).
   \end{equation}
   Using standard limiting arguments, we conclude that~\eqref{2.16p} holds for the actual solution to~\eqref{Sistema1FaseLinear}.

\textsc{Part~3:} We note that $\rho_{11}\Delta y_t\in L^2(Q)$ and, accordingly, we also have~$\rho_{17}\Delta y_t\in L^2(Q)$.
Furthermore,
   $$
\left|(\rho_{17}\Delta y_t)_t \right|\leq C\left(\rho_{17} |\Delta y_{tt}| + \rho_{11} |\Delta y_t| \right) \in L^2(Q).
   $$
This allows to deduce (again first for the Galerkin approximations and then for the solution to~\eqref{Sistema1FaseLinear}) that $\rho_{17}\Delta y_t \in L^\infty(0,T;L^2(\Om))$ and
   $$
\sup_{\left[0,T\right]}\left[\jnt_\Om \rho_{17}^2 |\Delta y_t|^2\,dx \right](t) \leq \jjntQ \left[ \rho_{17}^2|\Delta y_t|^2 + \left|\left(\rho_{17} \Delta y_{t}\right)_t\right|^2  \right] \leq S(y_0,v,f).
   $$

\textsc{Part~4:} Let us multiply (\ref{Sistema1FaseLinear}) by~$\rho_9^2\Delta^2 y$ and let us integrate in space.
   Then, after some manipulations, we see that
   $$
\begin{array}{l} \displaystyle
   \frac{1}{2}\frac{d}{dt}\left[\jnt_\Om \rho_9^2 |\Delta y|^2\,dx \right] + \sigma \jnt_\Om \rho_9^2 |\nabla \Delta y|^2\,dx
\\ \noalign{\smallskip} \displaystyle
\hspace{1.cm}
\leq C\bigg\{\jnt_\Om \left[\rho_9^2 |\nabla v|^2 + \rho_7^2|v|^2 \right] \,dx 
+ \jnt_\Om \left[\rho_9^2 |\nabla f|^2 + \rho_7^2 |f|^2 \right]dx
+\jnt_\Om \rho_7^2\left[y_t^2 + (\Delta y)^2 \right]dx \bigg\} 
\\ \noalign{\smallskip} \displaystyle
\hspace{1.cm} 
\ \ + \frac{\sigma}{2}\jnt_\Om \rho_9^2|\nabla \Delta y|^2\,dx .
\end{array}   
   $$
   In particular, we deduce the estimate
   $$
\jjntQ \rho_9^2 |\nabla \Delta y|^2\,dx \leq S(y_0,v,f).
   $$
   
   Let us apply the Laplace operator to both sides of~(\ref{Sistema1FaseLinear}).
   We get the identity $\Delta y_t -\sigma\Delta^2 y = \chi_\omega \Delta v + \Delta f$.
   By multiplying by $-\rho_{11}^2\Delta^2 y$ and proceeding as before, we see that
   $$
\begin{array}{l} \displaystyle
\frac{1}{2}\frac{d}{dt}\left[\jnt_\Om \rho_{11}^2 |\nabla \Delta y|^2\,dx \right] 
+ \sigma \jnt_\Om  \rho_{11}^2 |\Delta^2 y|^2\,dx 
\\ \noalign{\smallskip} \displaystyle
\hspace{1.cm} 
\ \leq C\bigg\{\jnt_\Om \rho_{11}^2|v|^2\,dx  + \jnt_\Om \left[\rho_{11}^2 |h|^2 +\rho_9^2 |\nabla \Delta y|^2 + \rho_{11}^2 |\Delta y_t|^2\right]\,dx \bigg\} 
+\ \frac{\sigma}{2}\jnt_\Om \rho_{11}^2 |\Delta^2 y|^2\,dx,
\end{array}
   $$
whence we find the estimate
   $$
\sup_{\left[0,T\right]}\left[ \jnt_\Om \rho_{11}^2 |\nabla \Delta y|^2\,dx \right](t) 
+ \jjntQ \rho_{11}^2 |\Delta^2 y|^2  \leq S(y_0,v,f)
   $$
for the Galerkin approximations.
   Once more, by standard limiting arguments, we see that this also holds for the solution to~\eqref{Sistema1FaseLinear}.

   Finally, in order to prove that $\rho_{13} \nabla \Delta y_t \in \left[L^2(Q)\right]^d$, it suffices to differentiate  (\ref{Sistema1FaseLinear}) with respect to time, multiply by $\rho_{13}^2 \Delta^2 y_t$ and work just as at the beginning of~\textsc{Part~4.}

   This ends the proof.

\section{Proof of Lemma~\ref{WellDefiniteness}}\label{Sec-App-B}

   We must prove that $H$ takes values in~$Z$.
   Obviously, it suffices to prove that the first component of~$H(y,v)$ belongs to~$F$.
   
   Let us put $H_1(y,v) = R(y,v) - B(y)$, with
   $$
R(y,v):= y_t - a(0)\Delta y -\chi_\omega v1_\omega \ \text{ and } \ B(y):= \nabla \cdot \left[ \left(a(\nabla y) - a(0)\right)\nabla y\right].
   $$
   From the definitions of~$Y$ and the norm of $F$, it is clear that $R(y,v) \in F$.
   
   Recall the notation $H_k = L^2(Q,\mu_k)$, where $d\mu_k = \rho_k^2 \,d(t,x)$.

\textsc{Claim~1:} $\Delta B(y) \in H_7$.
   
   It is straightforward to deduce that
   \begin{equation} \label{Pontwise1}
\begin{array}{l}
   |\Delta B(y)| \leq C\Big(|\nabla y |^r|D^4 y| +|\nabla y|^{r-2}|D^2y|^3 + |\nabla y|^{r-1}|D^2 y||D^3 y|   
   + |D^2 y||D^3 y| 
   \\ \noalign{\smallskip} \displaystyle
   \hspace{1.3cm}\,+ |\nabla y||D^2 y||D^3 y| + |\nabla y||D^4 y| + |D^2 y|^3 + |\nabla y||D^2 y|^3 \Big)
   \\ \noalign{\smallskip} \displaystyle
   \hspace{1.15cm}\ = C\left(B_1 +\cdots + B_8 \right) .
\end{array}   
   \end{equation}
   Taking into account the Sobolev embedding $H^2(\Om) \hookrightarrow L^\infty(\Om)$ and~\eqref{eq1}, we have
   \begin{equation} \label{B1}
\begin{array}{l} \displaystyle
   \jjntQ \rho_7^2 |B_1|^2  \leq \jnt_0^T e^{2s\alpha_2/\ell} \ell^7 \left(\jnt_\Om |\nabla y|^{2r}|D^4 y|^2 \,dx\right)\,dt
   \leq \jnt_0^T e^{2s\alpha_2/\ell} \ell^7 \|\nabla y\|_{L^\infty}^{2r}\|D^4 y\|^2 \,dt
   \\ \noalign{\smallskip} \displaystyle
   \phantom{\jjntQ \rho_7^2 |B_1|^2} \leq C \left[\sup_{\left[0,T\right]}\|\left( \rho_{13} y \right) (t )\|_{H^3(\Om)} \right]^{2r} 
   \int_0^T \| \left(  \rho_{11} y \right)(t)\|_{H^4(\Om)}^2  < +\infty
\end{array}   
   \end{equation}
and similar estimates can be found for the weighted integrals of~$B_2$ to~$B_8$.

\begin{remark}\label{rem3.6}
   The main reason why we have to impose $d \leqslant 3$ is that, in these estimates, we need the following results from Sobolev space theory:
\begin{itemize}
   \item[(a)] $H^1(\Om) \hookrightarrow L^6(\Om)$, with a dense and continuous embedding.
   \item[(b)] $H^2(\Om) \hookrightarrow L^\infty(\Om)$, with a dense and compact embedding.
\end{itemize}
   For $d \leq 3$, they are well known and can be easily deduced from the general results presented in~\cite{adams2003sobolev}.
\end{remark}

   Consequently, $\Delta B(y) \in H_7$ and \textsc{Claim~1} holds.

\textsc{Claim~2:} $ \nabla B(y)_t \in H_{11}$.

   We have the pointwise estimate
   \begin{equation} \label{Pointwise2}
\begin{array}{l}
   |\nabla \partial_t B(y)|\leq C\big(|\nabla y|^r |D^3 y_t| + |\nabla y|^{r-1}|\nabla y_t| |D^3 y| + |\nabla y|^{r-1}|D^2 y||D^2 y_t|
   + |\nabla y|^{r-2}|\nabla y_t||D^2 y|^2
   \\ \noalign{\smallskip} \displaystyle
   \hspace{2.2cm} \, + |\nabla y_t||D^2 y|^2 +  |\nabla y||\nabla y_t||D^2 y|^2 + |\nabla y||D^2 y||D^2 y_t| + |\nabla y||D^3 y_t| 
   + |\nabla y_t||D^3 y|
   \\ \noalign{\smallskip} \displaystyle
   \hspace{2.2cm} \, + |\nabla y||\nabla y_t||D^3 y| + |D^2 y||D^2 y_t| \big)
   \\ \noalign{\smallskip} \displaystyle
   \hspace{1.5cm} \,= C\left(C_1 + \cdots + C_{11} \right) .
\end{array}   
   \end{equation}
   
   For example,
   \begin{equation} \label{C1}
\begin{array}{l} \displaystyle
   \jjntQ \rho_{11}^2 |C_1|^2  \leq \jnt_0^T e^{2s\alpha_2/\ell} \ell^{11} 
   \left(\jnt_\Om |\nabla y|^{2r} |D^3 y_t|^2 \,dx\right) \,dt
   \leq C\jnt_0^T e^{2s\alpha_2/\ell} \ell^{11}\|\nabla y\|_{L^\infty}^{2r}\|D^3 y_t\|^2 \,dt
   \\ \noalign{\smallskip} \displaystyle
   \phantom{jjntQ \rho_{11}^2 |C_1|^2  } \leq C\left[ \sup_{\left[0,T\right]} \|(\rho_{13}y)(t)\|_{H^3(\Om)} \right]^{2r}
   \int_0^T \| (\rho_{13} y_t)(t)\|_{H^3(\Om)}^2  < +\infty
\end{array}   
   \end{equation}
and similar estimates can be deduced for~$C_2$ to~$C_{11}$.

   We easily find that~$\nabla B(y)_t \in H_{11}$ and, therefore, \textsc{Claim~2} also holds.

\textsc{Claim~3:} $B(y)_{tt} \in H_{17}$.

   In this case, we have 
   \begin{equation} \label{Pointwise3}
\begin{array}{l}
   |B(y)_{tt}| \leq C\big(|\nabla y|^{r}|D^2 y_{tt}| + |\nabla y|^{r-1}|\nabla y_t||D^2 y_t| + |\nabla y|^{r-1}|D^2 y||\nabla y_{tt}| 
   \\ \noalign{\smallskip} \displaystyle
   \hspace{2.2cm} \,+\, |\nabla y|^{r-2}|D^2 y||\nabla y_t|^2 + |\nabla y||D^2 y_{tt}| + |\nabla y_t||D^2 y_t| 
   \\ \noalign{\smallskip} \displaystyle
   \hspace{2.2cm} \,+\, |\nabla y||\nabla y_t||D^2 y_t| + |D^2 y||\nabla y_{tt}| + |\nabla y||D^2 y||\nabla y_{tt}| 
   \\ \noalign{\smallskip} \displaystyle
   \hspace{2.2cm} \,+\, |D^2 y||\nabla y_t|^2 + |\nabla y||D^2 y||\nabla y_t|^2 \big)
   \\ \noalign{\smallskip} \displaystyle
   \hspace{1.2cm} \ = C(D_1 + \cdots + D_{11}).
\end{array}   
   \end{equation}
   
   As before, it suffices to prove that the $D_j$ belong $H_{17}$ for $j = 1, \dots ,4$.
   For instance, in the case of~$D_1$, we have:
   \begin{equation} \label{D1}
   \begin{array}{l} \displaystyle
   \jjntQ \rho_{17}^2 |D_1|^2  \leq C \jnt_0^T 
   e^{2s\alpha_2/\ell} \ell^{17}\|\nabla y\|_{L^\infty}^{2r}\|\Delta y_{tt}\|^2 \,dt
   \\ \noalign{\smallskip} \displaystyle
   \phantom{\jjntQ \rho_{17}^2 |D_1|^2  } \leq C \left[\sup_{\left[0,T\right]}\|(\rho_{13} y)(t)\|_{H^3(\Om)} \right]^{2r}
   \jjntQ \rho_{17}^2 (\Delta y_{tt})^2 < +\infty .
\end{array}   
   \end{equation}
   
   Similar estimates can be obtained for the other $D_j$.

   The remaining terms in the norm of $B(y)$ in~$F$ are even easier to estimate.
   Thus, we get that~$B(y) \in F$ and~$H : Y \mapsto Z$ is well-defined.
   
   The continuity of~$H$ is deduced arguing in a similar way.
   For brevity, the details are left to the reader.

\section{Proof of Lemma~\ref{ContinuousAndOnto}}\label{Sec-App-C}

   The fact that $\Lambda$ is continuous is immediate from the definition of the norms of $Y$ and~$Z$. 
   Next, we address the surjectiveness of $\Lambda$.
   
   Let us fix $(f,y_0) \in Z$.
   By Theorem~\ref{ControlOfLinearProblem}, there exists a state-control pair $(y,v)$ corresponding to the data $f$ and~$y_0$, with $\sigma$ replaced by~$a(0)$.
   From the estimates in~Corollary \ref{RegularityForLinear1PhaseProblem}, we have~$(y,v) \in Y$.
   Furthermore,
   \begin{equation} \label{surjectiveLambda}
\Lambda (y,v) = \left( y_t - a(0)\Delta y -\chi_\omega v, y(0\,,\cdot) \right) = (f,y_0).
   \end{equation}
   This shows that $(f,y_0)$ belongs to the rank of~$H$, that is, $\Lambda$ is surjective.
   
   Moreover, we also deduce from Corollary~\ref{RegularityForLinear1PhaseProblem} that there exists $M>0$ with
   $$
\|(y,v)\|_Y \leq M^{-1}\|(f,y_0)\|_Z = M^{-1}\|\Lambda(y,v)\|_Z.
   $$
   After the change of variable in time~$t^\prime=t\sigma$, we can rewrite the Carleman inequality~\eqref{CI-new} in the form
   \begin{equation} \label{CI-new2}
   I(s,\lambda;\varphi) \leq \overline{C}\left(\jjntQ e^{-2s\alpha}|F|^2\  + \jjntomsT e^{-2s\alpha}(s\zeta)^3|\varphi|^2 \right),
   \end{equation}
   where $\overline{C}=C(\Om, \omega_0, T)\max \{\sigma, 1/\sigma\}\left(\min \{\sigma,1/\sigma\}\right)^{-1}$.
   Since 
   \begin{equation} \label{NormOfY_new}
\|(y,v)\|_Y^2 := \|v\|_3^2 + \|v_t\|_7^2 + \|\Delta v\|_7^2 + \|D^4 v\|_{15}^2 + \|\Delta v_t\|_{15}^2 + \|v_{tt}\|_{15}^2 +\|\Lambda(y,v)\|_Z^2 ,
   \end{equation}
using the Carleman inequality constant $\overline{C}$ in~Theorem~\ref{ControlOfLinearProblem} and~Corollary~\ref{RegularityForLinear1PhaseProblem}~(b) corresponding to~$\sigma=a(0)$, we have from \eqref{NormOfY_new} that
   $$
\|(y,v)\|_Y \leq \overline{C} \left(\max \left\{ \frac{1}{a(0)},1 \right\} \right) \|\Lambda(y,v)\|_Z + \|\Lambda(y,v)\|_Z.
   $$
   Consequently, $\|\Lambda(y,v)\|_Z \geq M\|(y,v)\|_Y$ with
   \[
M= \left[C(\Om, \omega_0, T)\max \{a(0), \frac{1}{a(0)}\}\left(\min \{a(0),\frac{1}{a(0)}\}\right)^{-1}\max\{\frac{1}{a(0)},1\}+1\right]^{-1} .
   \]
   
   This ends the proof.
   
\section{Proof of Lemma~\ref{StrictDifferentiabilityAtTHeOrigin}}\label{Sec-App-D}

  Without loss of generality, we can assume that $r \geq 4$.
   Let us take $(y,v), (\overline{y},\overline{v}) \in Y$ and let~$0< \delta \leq 1$ be given.
   We also assume that $\|(y,v)\|_Y \leq \delta$ and~$\|(\overline{y},\overline{v})\|_Y \leq \delta$.
   
   Let us put
   \begin{equation} \label{strictDiffStep1}
   H(y,v) - H(\overline{y},\overline{v}) - \Lambda(y-\overline{y}, v- \overline{v}) = \left( -\Phi(y,\overline{y}),0 \right),
   \end{equation}
where
   \begin{equation} \label{definitionOfPhi}
   \Phi(y,\overline{y}) := \nabla \cdot \left[\left( a(\nabla y) - a(0) \right) \nabla \left(y-\overline{y}\right) + \left( a(\nabla y) - a(\nabla \overline{y})\right) \nabla \overline{y} \right].
   \end{equation}
   Then,
   \begin{equation} \label{PhiConnectedWithH}
   \|H(y,v) - H(\overline{y},\overline{v}) - \Lambda(y-\overline{y},v- \overline{v})\|_Z = \|\Phi(y,\overline{y})\|_F.
   \end{equation}
   We will prove that, for some $K > 0$,
   \begin{equation} \label{NormFofPhi}
   \|\Phi(y,\overline{y})\|_F \leq K \delta \|(y-\overline{y},v-\overline{v})\|_Y.
   \end{equation}
   In view of (\ref{PhiConnectedWithH}), this will suffice.
   
   The first step is to show that
   \begin{equation} \label{FirstStepPhi}
   \|\Delta \Phi(y,\overline{y})\|_7 \leq B\delta \|(y-\overline{y},v-\overline{v})\|_Y ;
   \end{equation}
   recall the definition of~$\|\cdot\|_7$ in~\eqref{3.1a}.
To begin with, we observe that
   $$
   \begin{array}{l}
      \Delta \Phi(y,\overline{y}) = \Delta\left[a(\nabla y) \right]\Delta\left(y-\overline{y}\right) + \left[a(\nabla y) - a(0) \right]\Delta^2(y-\overline{y})  
      \\ \noalign{\smallskip} \displaystyle
      \hspace{1.6cm} \,+\, 3\nabla\left[a(\nabla y) \right]\nabla \Delta (y-\overline{y}) + \nabla \Delta \left[ a(\nabla y) \right] \nabla (y-\overline{y}) 
      \\  \noalign{\smallskip} \displaystyle
      \hspace{1.6cm} \,+\,2\tr\left(D^2\left[a(\nabla y) \right]^{T}D^2( y -\overline{y})\right) + \Delta\left[a(\nabla y) - a(\nabla\overline{y}) \right]\Delta \overline{y} 
      \\ \noalign{\smallskip} \displaystyle
      \hspace{1.6cm} \,+\, \left[a(\nabla y) - a(\nabla\overline{y})\right]\Delta^2\overline{y} +3\nabla\left[a(\nabla y) - a(\nabla\overline{y})\right]\nabla\Delta \overline{y}
      \\ \noalign{\smallskip} \displaystyle
      \hspace{1.6cm} \,+\, \nabla \Delta\left[a(\nabla y) - a(\nabla\overline{y})\right] \nabla\overline{y} + 2\tr\left( D^2\left[a(\nabla y) - a(\nabla\overline{y}) \right]^{T}D^2\overline{y} \right) ,
   \end{array}
   $$
which yields
   $$
   \begin{array}{l}
      |\Delta \Phi (y,\overline{y}) | \leq C\big\{|D^2\left[ a(\nabla y) \right]||D^2(y-\overline{y})| + |a(\nabla y) - a(0)||\Delta^2 (y - \overline{y})| 
      \\ \noalign{\smallskip} \displaystyle
      \hspace{1.9cm} \,+\, |\nabla\left[ a(\nabla y)\right]||\nabla\Delta (y-\overline{y})| + |\nabla \Delta \left[ a(\nabla y) \right]||\nabla (y-\overline{y})|  
      \\ \noalign{\smallskip} \displaystyle
      \hspace{1.9cm} \,+\, |D^2\left[ a(\nabla y) -a(\nabla \overline{y}) \right]||D^2\overline{y}| + |a(\nabla y) - a(\nabla \overline{y})||\Delta^2 \overline{y}| 
      \\ \noalign{\smallskip} \displaystyle
      \hspace{1.9cm} \,+\, |\nabla\left[a(\nabla y) - a(\nabla\overline{y}) \right]||\nabla\Delta \overline{y}| + |\nabla\Delta \left[ a(\nabla y) - a(\nabla \overline{y})\right]||\nabla \overline{y}| \big\}
      \\ \noalign{\smallskip} \displaystyle
      \hspace{1.6cm} \,=\, C\left(E_1 + \cdots + E_{8} \right). 
   \end{array}
   $$

   Using hypothesis~\textbf{H2} and then the estimates in~Corollary~\ref{RegularityForLinear1PhaseProblem}, we can deal with all the terms $E_1,...,E_8$.
   Thus, let us show what can be done with $E_5$;
   the other $E_j$ can be bounded likewise.
   
   By adding and subtracting adequate terms and using the Mean Value Theorem, we get:
   $$
   \begin{array}{l}
   E_5 \leq C\big[ \left(|\nabla y|^{r-3} + |\nabla \overline{y}|^{r-3} \right)|D^2 y|^2|D^2 \overline{y}||\nabla (y-\overline{y})|   
   \\ \noalign{\smallskip} \displaystyle
   \hspace{1.2cm} \,+\, \left(|\nabla y|^{r-2} + |\nabla \overline{y}|^{r-2} \right)|D^3 y||D^2 \overline{y}||\nabla (y-\overline{y})|  
   \\ \noalign{\smallskip} \displaystyle
   \hspace{1.2cm} \,+\, |\nabla\overline{y}|^{r-2}|D^2 y||D^2 \overline{y}||D^2(y-\overline{y})| 
   +  |\nabla \overline{y}|^{r-1}|D^2 \overline{y}||D^3(y-\overline{y})|  
   \\ \noalign{\smallskip} \displaystyle
   \hspace{1.2cm} \,+\,  |\nabla\overline{y}|^{r-2}|D^2\overline{y}|^2|D^2(y-\overline{y})| 
   + |D^2 y|^2|D^2 \overline{y}||\nabla (y-\overline{y})| \\ \noalign{\smallskip} \displaystyle
   \hspace{1.2cm} \,+\, |D^3 y||D^2 \overline{y}||\nabla (y-\overline{y})| + |D^2 y||D^2 \overline{y}||D^2(y-\overline{y})| 
   \\ \noalign{\smallskip} \displaystyle
   \hspace{1.2cm} \,+\,  |D^2 \overline{y}||D^3(y-\overline{y})| + |D^2\overline{y}|^2|D^2(y-\overline{y})| \big]
   \\ \noalign{\smallskip} \displaystyle
   \hspace{0.52cm} = C(E_{5,1} + \cdots + E_{5,10}) .
   \end{array}
   $$
   For example,
   \begin{equation} \label{E51}
   \begin{array}{l} \displaystyle
      \jjntQ \rho_7^2 |E_{5,1}|^2 \leq C\jnt_0^T e^{2s\alpha_2/\ell} \ell^{7}
      \left( \|\nabla y\|_{L^\infty(\Om)}^{2(r-3)} + \|\nabla \overline{y}\|_{L^\infty(\Om)}^{2(r-3)} \right)
      \| \nabla(y-\overline{y})\|_{L^\infty}^2\|D^2 y\|_{L^6}^4 \|D^2 \overline{y}\|_{L^6}^2 \,dt
      \\ \noalign{\smallskip} \displaystyle
      \phantom{\jjntQ \rho_7^2 |E_{5,1}|^2 } \leq C\left[ \sup_{\left[0,T\right]} \left( \|\rho_{13} y\|_{H^3(\Om)} 
      + \|\rho_{13} \overline{y}\|_{H^3(\Om)} \right)(t) \right]^{2r} 
      \left[ \sup_{\left[0,T\right]} \| \rho_{13} (y-\overline{y})\|_{H^3(\Om)}(t) \right]^2
      \\ \noalign{\smallskip} \displaystyle
      \phantom{\jjntQ \rho_7^2 |E_{5,1}|^2 } \leq C\delta^{2r}\|(y-\overline{y},v-\overline{v})\|_Y^2
   \end{array}
   \end{equation}
and we can obtain very similar inequalities for the other $E_{5,i}$.

   As a final consequence, we get~\eqref{FirstStepPhi} and the proof is done.

\section*{Acknowledgements}
   EFC was partially financed by MINECO (Spain), Grant~PID2020-114976GB-I00.
   DM and YT were financed in part by Coordena\c{c}\~ao de Aperfei\c{c}oamento de Pessoal de N\'ivel Superior - Brasil (CAPES) - Finance code 001.


\end{document}